\documentclass[11pt]{article}
\usepackage[applemac]{inputenc}
\usepackage{microtype}
\usepackage[french,english]{babel}
\usepackage{amssymb}
\usepackage{amsmath}
\usepackage{amsthm}
\usepackage{graphicx}
\usepackage{libertine}
\usepackage{newtxmath}
\usepackage{xcolor}
\definecolor{myblue}{rgb}{0.153,0.322,0.706}
\usepackage[colorlinks,linkcolor=myblue,urlcolor=myblue,citecolor=myblue,breaklinks=true]{hyperref}
\usepackage{paracol}
\usepackage{endnotes}
\usepackage{geometry}
\usepackage{multicol}
\usepackage{enumitem}

\newcommand{\be}{\begin{equation}}
\newcommand{\ee}{\end{equation}}
\newcommand{\ra}{\rightarrow}
\newcommand{\bV}{\bar{V}}
\newcommand{\bZ}{\bar{Z}}
\newcommand{\bm}{\bar{m}}
\newcommand{\bsigma}{\bar{\sigma}}
\newcommand{\bF}{\bar{F}}
\newcommand{\bW}{\bar{W}}
\newcommand{\bnu}{\bar{\nu}}
\newcommand{\bw}{\bar{w}}
\newcommand{\bS}{\bar{S}}

\colorlet{mylightgray}{gray!40}
\colseprulecolor{mylightgray}
\normalcolseprulecolor[1]
\makeatletter
\def\enoteformat{%
  \rightskip\z@ \leftskip\z@ \parindent=1.8em
  \leavevmode{\setbox\z@=\lastbox}\llap{[\![\theenmark]\!]\enskip}%
}
\makeatother

\newtheoremstyle{myplain}
{5pt}			
{5pt}			
{\em}			
{\parindent}	
{\scshape}		
{.}			
{.5em}		
{\thmname{#1}\thmnumber{ #2}\thmnote{~{(#3)}}}

\theoremstyle{myplain}
\newtheorem{conditionth}{Condition}
\newtheorem{theoremen}{Theorem}
\newtheorem{theoremfr}{Th\'eor\`eme}

\begin{document}

\newgeometry{margin=1.2in,top=1in,bottom=1in}

\title{\textbf{On a new limit theorem in probability theory}\\
\mbox{(Sur un nouveau th\'eor\`eme-limite de la th\'eorie des probabilit\'es)}
}
\author{Harald Cram\'er (1893-1985)\\
\textit{Stockholm, Sweden}\\
\ \\
\textit{Translated by}\\ 
Hugo Touchette\\
\emph{National Institute for Theoretical Physics \emph{(}NITheP\emph{)}, Stellenbosch, South Africa}
}
\date{15 March 2018}
\maketitle

\noindent \emph{Original article}: H. Cram\'er, Sur un nouveau th\'eor\`eme-limite de la th\'eorie des probabilit\'es, Colloque consacré à la théorie des probabilités, Actualit\'es scientifiques et industrielles 736, 2-23, Hermann \& Cie, Paris, 1938.

\bigskip

\noindent\emph{Reprinted in}: H. Cram\'er, \href{http://www.springer.com/gp/book/9783642396847}{\emph{Collected Works}}, A. Martin-L\"of (Ed.), Vol. II, Springer, Berlin, 1994, p.~895-913.

\vspace*{0.2in}

\begin{paracol}{2}
\centering
\includegraphics[width=0.44\textwidth,angle=0.7]{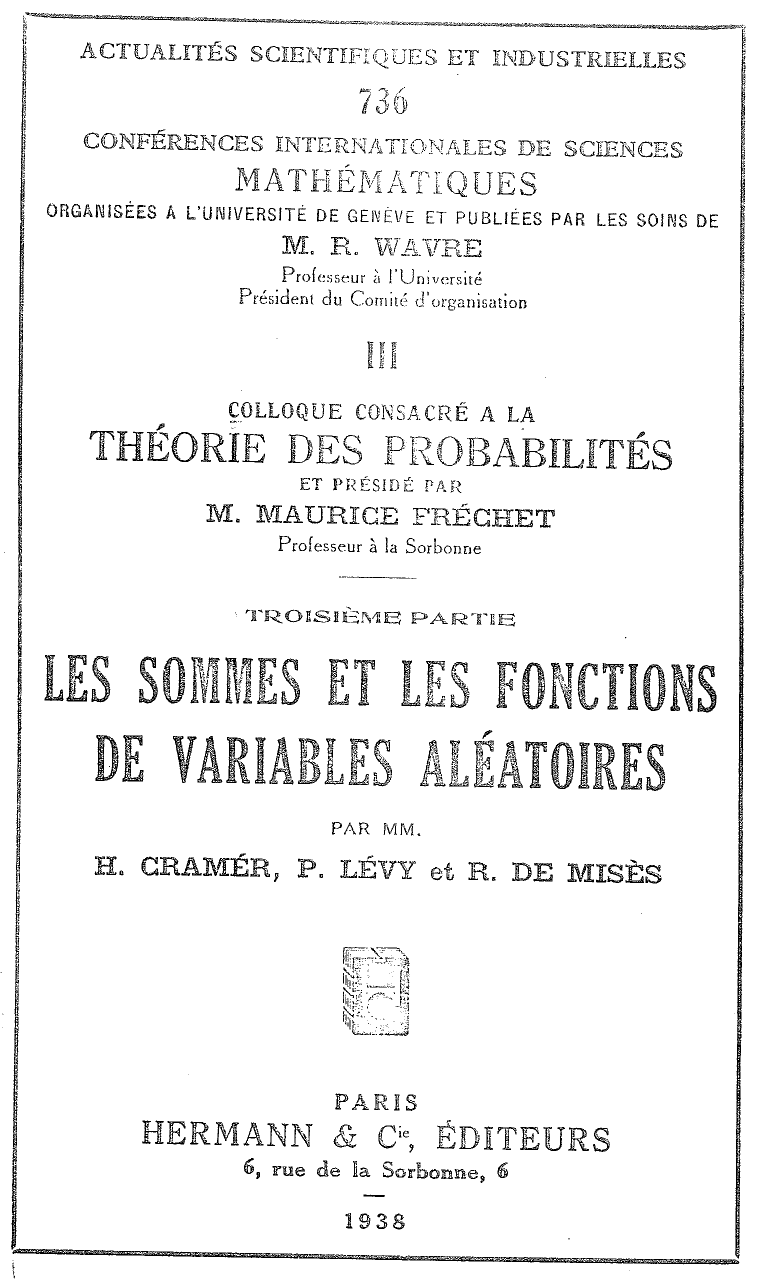}

\switchcolumn
\includegraphics[width=0.41\textwidth,angle=0.7]{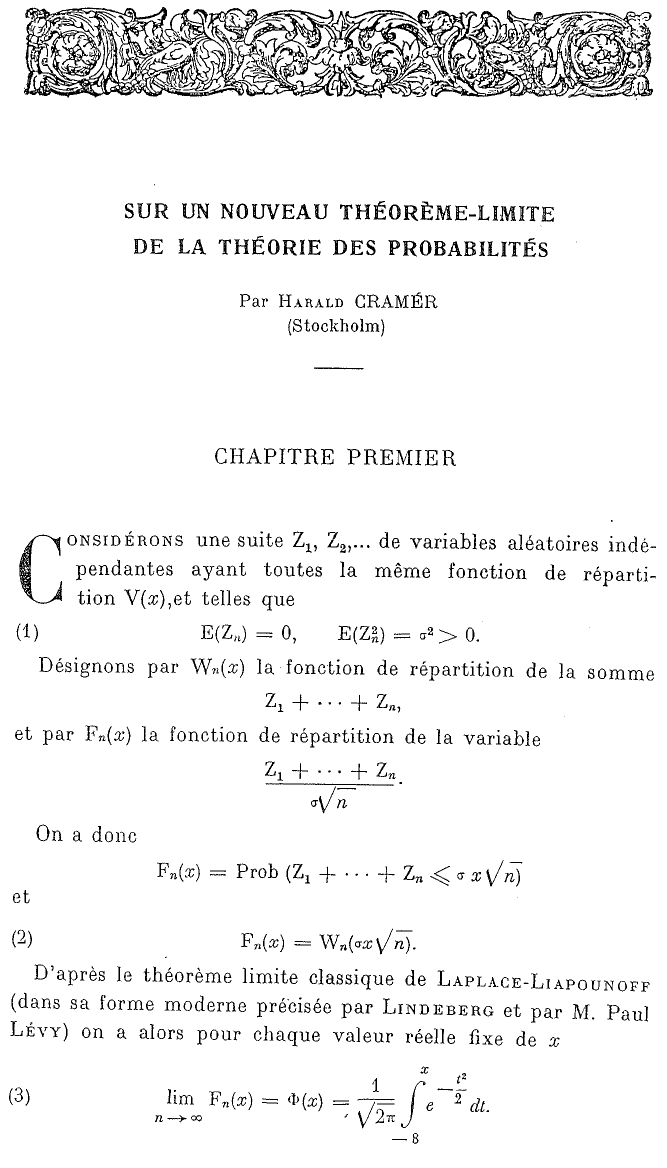}
\end{paracol}

\newpage

\newgeometry{margin=0.75in,top=1in,bottom=1in}
\fontsize{10pt}{12pt}\selectfont

\setlength{\columnsep}{0.3in}
\begin{multicols}{2}
\setlength{\columnseprule}{0pt}

\section*{Introduction}

The following is a translation of Harald Cram\'er's article, ``On a new limit theorem in probability theory'', published in French in 1938 \cite{cramer1938} and deriving what is considered by mathematicians to be the first large deviation result (Theorem 6 in the text). There is no electronic version of Cram\'er's article available from Hermann Editions, the publisher of the proceedings in which Cram\'er's contribution appears, although \href{http://books.google.com/books/about/Les_sommes_et_les_fonctions_de_variables.html?id=1jIgAAAAIAAJ&redir_esc=y}{Google Books has a partially-readable version}. My hope is that this translation will help disseminate this historically important work, 80 years after its publication. The conference held in Geneva in 1937 for which it was written is briefly mentioned by Cram\'er in his ``Personal Recollections'' \cite{cramer1976} (see Sec.~4.9, p.~528).

I have been actively looking for an English translation of this article, but could not find any, even with help from researchers familiar with it. The article is reprinted in volume II of Cram\'er's \emph{Collected Works} \cite{cramer1994}, but in the original French. Moreover, there happens to be a Russian version in \href{http://www.mathnet.ru/eng/umn}{\emph{Uspekhi Matematicheskikh Nauk}} \cite{cramer1944}, translated in 1944 by none other than Gnedenko, but it is not translated in English in the \href{http://iopscience.iop.org/journal/0036-0279}{\emph{Russian Mathematical Surveys}}, the latter having started only in 1960. Finally, the bibliographies and obituaries of Cram\'er I could find \cite{blom1984,blom1987,leadbetter1988} list only the French 1938 article. Needless to say, I would be grateful if someone could point an English translation to me.

\section*{Context}

The importance of Cram\'er's article has been discussed before, notably by Le Cam \cite[Sec.~4]{cam1987}, so I decided not to provide a commentary myself. His results are also now standard in large deviation theory, and are well explained, for example, in the textbooks of Dembo and Zeitouni \cite{dembo1998} and of den Hollander \cite{hollander2000}, in addition to several reviews \cite{ellis1999,amann1999,varadhan2008,touchette2009}.

I have added some notes in the English text (with the double brackets $[\![\ ]\!]$) only to highlight some important parts and results, and to comment on minor translation choices I have made. A brief synopsis of the article and some sources about Cram\'er's life and work are also included in the next sections. My feeling is that the importance of Cramer's article speaks, in any case, for itself. Like many classics, it should be read in its original version, in the voice of the period in which it was written.

I will only say, to put some historical context, that Cram\'er's results were definitively novel (Blom \cite{blom1987} calls them ``pioneering'') and somewhat out of the fashion of the time in probability theory for trying to go beyond the central limit theorem \cite{cramer1976,cam1986}. The conference in Geneva, presided by Fr\'echet, where these results were presented is also considered as a turning point or even an ``epiphany''\footnote{Laurent Mazliak, private communication, 2018.} in the history of probability theory, as a more rigorous approach was being adopted to study more abstract problems \cite{bustamente2015}.  This probably explains altogether why Cram\'er goes into such length in his first two chapters to explain the new large deviation problem that he wants to solve, how it relates to the central limit theorem, and why new methods going beyond that theorem are needed for tackling the problem.

It is clear that Cram\'er's inspiration for this problem came from his own work on the central limit theorem \cite{cramer1976,cam1986,cam1987} and, more importantly, from his experience of insurance and actuarial mathematics, on which he worked all his life \cite{cramer1976,martin2014}. The problem, as he explains in his recollections \cite[p.~514]{cramer1976}, is that
\begin{quotation}
\noindent for the risk problems in which I was interested, it was not enough to know that a certain probability distribution [for a sum of random variables] was approximately normal; it was necessary to have an idea of the magnitude of the error involved in replacing the distribution under consideration by the normal one. 
\end{quotation}
In other words, Cram\'er was looking for ``corrections'' to the central limit theorem, which is what he presents in his 1938 article. He considers a sum of $n$ independent and identically distributed random variables, representing for an insurance company a sum of $n$ claims, and tries to obtain an estimate of the probability distribution for ``large'' sum values that scale with $n$, as opposed to $\sqrt{n}$, which is the ``normal'' scaling considered in the central limit theorem. This is the large deviation problem, which he solves by providing (in Theorem 6) an asymptotic expansion of the distribution valid as $n\ra\infty$, whose dominant term is a decaying exponential in $n$.

The importance of this simple result is now clear, 80 years after its publication, and deserves indeed to be called ``pioneering'', having led to the development of a whole theory of large deviations, which finds applications in fields as diverse as queueing theory, information theory, statistical estimation, hypothesis testing, and statistical physics \cite{dembo1998,hollander2000,touchette2009}. Many of the techniques used by Cram\'er in 1938 are also now standard in probability theory. They include, in particular, the Esscher transform \cite{esscher1932}, a change of measure widely used in large deviation theory, statistics and rare event simulations under the name of exponential family, exponential tilting or Cram\'er transform. The dominant exponential term appearing in the asymptotic expansion found by Cram\'er is also the cornerstone of large deviation theory, referred to as the large deviation principle. The exponent appearing in that term is called the rate function or the Cram\'er function.

The importance of Cram\'er's article can also be measured, in the end, by the time it took for researchers to continue what Cram\'er had started. Results on large deviations began to appear after 1938 only around the 1960s, with works from Linnik, Petrov, and Sanov, to mention a few (see \cite{cam1987} for references). These results, however, do not form any general theory, but concern themselves with different sums of independent random variables or different conditions for results similar to Cram\'er's to apply. The starting point of the theory of large deviations, as we know it today, is considered to be the work of Donsker and Varadhan (see \cite{varadhan2008} for references) published in the 1970s on large deviations of Markov processes. From that time, and especially from the 1980s, the theory has grown to become one of the most active subjects in probability theory and in statistical physics.

In physics, Boltzmann is credited for having derived in 1877 the first large deviation result \cite{ellis1999}. However, there is no doubt that it is Cram\'er who initiated the mathematical study of large deviations.

		
\section*{Acknowledgments}

My thanks go to Arnaud Guyader for encouragements during the project and for carefully proofreading the translation. All remaining errors are, of course, mine. I also thank Richard S.\ Ellis for giving me a first copy of Cram\'er's article, Laurent Mazliak for useful comments and references about the 1937 Geneva conference, and Rapha\"el Chetrite for offering to participate in the project and for his encouragements.

\section*{Copyright}

I wrote to \'Editions Hermann in Paris to get permission for this translation, to post it on the arXiv and, ultimately, to get permission to publish it somewhere. I first sent them an email in 2012 and then wrote a letter in 2013, but did not receive an answer to either. 

I would be happy to be contacted by the publisher, should they see the translation and have any qualms with it. My understanding is that, although the paper is 80 years old, only 33 years have elapsed since Cram\'er's death in 1985, which is below the 50 years needed for his work to go into the public domain. I do not claim copyright for the translation and nor does the arXiv under its normal license.

\section*{Synopsis}

\begin{itemize}[leftmargin=*]
\item Chapter I:
\begin{itemize}[leftmargin=*]
\item Problem definition: Find the distribution of a sum of IID random variables with common cumulative distribution function (CDF) $V(x)$.
\item Explains the central limit theorem (CLT) after Eq.~(2).
\item Defines the large deviation problem in the paragraphs coming after Eq.~(4) up to Condition A.
\item Discusses whether the CLT can be used to solve this problem; the answer is no.
\item Defines Condition A, requiring the existence of the generating function of the CDF $V(x)$. This condition is fundamental in large deviation theory and is now referred to as Cram\'er's condition.
\end{itemize}

\item Chapter II:~%
\begin{itemize}[leftmargin=*]
\item Introduces the change of CDF, referred to as the Esscher transform, underlying the results to be discussed after.
\item Relates the original CDF $V(x)$ and transformed CDF $\bar V(x)$ (relations between the `non-bar' and `bar' quantities).
\item Relates, in particular, the characteristic functions of the original and transformed CDFs; see Eqs.~(9) and (10).
\end{itemize}

\item Chapter III:
\begin{itemize}[leftmargin=*]
\item Studies the cumulants and give CLT-type error estimates.
\item Theorem 1: Error estimate for the CLT, based on the scaling $x\sim \sqrt{n}/\log n$.
\item Theorem 2: Similar error estimate for $x\sim n^{1/6}$.
\item Theorem 3: Rewriting of Theorem 2, essentially.
\end{itemize}

\item Chapter IV:
\begin{itemize}[leftmargin=*]
\item Brief discussion of the binomial case. Not essential.
\item Discusses relations with other results of Khintchine, Levy, etc.
\end{itemize}

\item Chapter V (main chapter):
\begin{itemize}[leftmargin=*]
\item Defines Condition B needed to get rid of the $\log n$ term in the CLT error estimates. This condition seems unnecessary now in view of more refined results published after 1938; see \cite[Sec.~4]{cam1987}.
\item Theorem 5: Refinement of Theorems 1-3 in which the $\log$ term is omitted to consider the scaling $x\sim\sqrt{n}$.
\item Theorem 6: Large deviation result for the scaling $x\sim n$, coming from the $\sqrt{n}$ term in $F_n(x)$. The ``rate exponent'' $\alpha$ is the Legendre transform of the cumulant function; see Eqs.~(27) and (29). This is what is referred now to as Cram\'er's Theorem. Condition B does not appear in the modern form of that theorem.
\end{itemize}

\item Chapter VI:
\begin{itemize}[leftmargin=*]
\item Studies continuous-time (process) generalizations.
\end{itemize}

\end{itemize}

\section*{More information}

\begin{itemize}

\item Biography: see \href{https://en.wikipedia.org/wiki/Harald_Cramer}{Wikipedia} and references therein.

\item Personal recollections of Cram\'er around the history of probability theory (from about 1920 to 1970), the Swedish school of probability theory, actuarial mathematics, and his own work: \cite{cramer1976}.

\item Similar recollections compiled by Wegman: \cite{wegman1986}.

\item Complete bibliography of Cram\'er: \cite{blom1984}.

\item Collected works of Cram\'er: \cite{cramer1994}.

\item Obituary notices: \cite{blom1987,leadbetter1988}.

\item Pictures of Cram\'er from the \href{http://dynkincollection.library.cornell.edu/biographies/823}{Dynkin library}.

\item A short audio clip from \href{http://hdl.handle.net/1813/17270}{Dynkin's audio collection}.
	
\end{itemize}

\bibliography{masterbib}
\bibliographystyle{myplain2-doi}
\end{multicols}

\newgeometry{margin=0.6in,top=0.9in,bottom=1in}
\fontsize{10pt}{12pt}\selectfont

\begin{paracol}{2}

\section*{Chapitre premier}

\switchcolumn

\section*{First chapter}

\switchcolumn*

Considérons une suite $Z_1,Z_2,\ldots$ de variables aléatoires indépen\-dantes ayant toutes la même fonction de répartition $V(x)$, et telles que 
\be
E(Z_n)=0,\qquad E(Z_n^2)=\sigma^2>0.
\label{eqfr1}
\ee

\switchcolumn

Consider a sequence $Z_1,Z_2,\ldots$ of independent random variables having the same cumulative distribution function \endnote{``Repartition function'' is translated as ``cumulative distribution function'' (CDF).} $V(x)$ and such that
\be
E(Z_n)=0,\qquad E(Z_n^2)=\sigma^2>0.
\label{eqen1}
\ee

\switchcolumn*

Désignons par $W_n(x)$ la fonction de répartition de la somme
\[
Z_1+\cdots+Z_n,
\]
et par $F_n(x)$ la fonction de répartition de la variable
\[
\frac{Z_1+\cdots+Z_n}{\sigma \sqrt{n}}.
\]

\switchcolumn

Denote by $W_n(x)$ the cumulative distribution function of the sum
\[
Z_1+\cdots +Z_n,
\]
and by $F_n(x)$ the cumulative distribution function of the variable
\[
\frac{Z_1+\cdots+Z_n}{\sigma \sqrt{n}}.
\]

\switchcolumn*

On a donc
\[
F_n(x)=\textrm{Prob}(Z_1+\cdots+Z_n\leq \sigma x\sqrt{n})
\]
et
\be
F_n(x)=W_n(\sigma x\sqrt{n}).
\label{eqfr2}
\ee

\switchcolumn

We \endnote{The third person singular ``on'' in French is translated as ``we'' throughout to avoid the correct but heavier ``one''.} thus have
\[
F_n(x)=\textrm{Prob}(Z_1+\cdots+Z_n\leq \sigma x\sqrt{n})
\]
and
\be
F_n(x)=W_n(\sigma x\sqrt{n}).
\label{eqen2}
\ee

\switchcolumn*

D'après le théorème limite classique de Laplace-Liapounoff (dans sa forme moderne précisée par Lindeberg et par M.\ Paul Lévy) on a alors pour chaque valeur réelle fixe de $x$
\be
\lim_{n\ra\infty} F_n(x)=\Phi(x)=\frac{1}{\sqrt{2\pi}}\int_{-\infty}^x e^{-\frac{t^2}{2}}\, dt.
\label{eqfr3}
\ee

\switchcolumn

Following the classical limit theorem of Laplace-Lyapunov \endnote{``Liapunov'' is used instead of ``Liapounoff''.} (in its modern version specified by Lindeberg and by Paul L\'evy) we thus have for each real value $x$
\be
\lim_{n\ra\infty} F_n(x)=\Phi(x)=\frac{1}{\sqrt{2\pi}}\int_{-\infty}^x e^{-\frac{t^2}{2}}dt.
\label{eqen3}
\ee

\switchcolumn*

Par ce théorème, on a donc une expression asymptotique (pour $n\ra\infty$) de la probabilité $F_n(x)$ de l'inégalité
\[
Z_1+\cdots+Z_n\leq \sigma x\sqrt{n}
\]
ou, ce qui revient au même, de la probabilité $1-F_n(x)$ de l'inégalité
\[
Z_1+\cdots+Z_n> \sigma x \sqrt{n}
\]
$x$ étant toujours un nombre réel \emph{indépendant de $n$}.

\switchcolumn

From this theorem \endnote{This is the central limit theorem.}, we thus have an asymptotic expression (for $n\ra\infty$) for the probability $F_n(x)$ of the inequality
\[
Z_1+\cdots+Z_n\leq \sigma x\sqrt{n}
\]
or, which amounts to the same, for the probability $1-F_n(x)$ of the inequality
\[
Z_1+\cdots+Z_n> \sigma x \sqrt{n}
\]
$x$ being as before a real number \emph{independent of $n$}.

\switchcolumn*

Il est alors naturel de se demander ce que deviennent ces probabilités \emph{lorsque $x$ peut varier avec $n$, en tendant vers $+\infty$ ou vers $-\infty$ quand $n$ croît indéfiniment}.

\switchcolumn

It is thus natural to ask what happens of these probabilities \textit{when $x$ can vary with $n$, going to $+\infty$ or to $-\infty$ when $n$ grows indefinitely}.

\switchcolumn*

Dans ces conditions, la relation \eqref{eqfr3} ne donne que le résultat évident
\[
\lim_{n\ra\infty} F_n(x)=
\left\{
\begin{array}{lcl}
1 & \text{quand} & x\ra +\infty,\\
0 & \text{''}& x\ra -\infty,
\end{array}
\right.
\]
qui exprime seulement que $F_n(x)$ tend vers les mêmes limites que $\Phi(x)$ lorsque $\ra\pm\infty$.

\switchcolumn

In these conditions, Relation \eqref{eqen3} only gives the evident result
\[
\lim_{n\ra\infty} F_n(x)=
\left\{
\begin{array}{lcl}
1 & \text{when} & x\ra +\infty,\\
0 & \text{''} & x\ra -\infty,
\end{array}
\right.
\]
which expresses only that $F_n(x)$ converges to the same limits as $\Phi(x)$ when $\ra\pm\infty$.

\switchcolumn*

Pour savoir si l'équivalence asymptotique de $F_n(x)$ et $\Phi(x)$ subsiste dans les conditions indiquées, on pourrait se proposer d'étudier les rapports
\be
\frac{1-F_n(x)}{1-\Phi(x)}\quad\text{pour}\quad x\ra +\infty,\tag{4a}
\label{eqfr4a}
\ee 
et
\be
\frac{F_n(x)}{\Phi(x)}\quad\text{pour}\quad x\ra -\infty.\tag{4b}
\label{eqfr4b}
\ee
Si $x$ est indépendant de $n$, il suit de \eqref{eqfr3} que ces rapports tendent tous les deux vers l'unité lorsque $n$ tend vers l'infini\,; il s'agit maintenant de savoir ce qui arrive quand $x$ devient infini avec $n$.

\switchcolumn

To see whether the asymptotic equivalence of $F_n(x)$ and $\Phi(x)$ remains under the mentioned conditions indicated, we could propose to study the ratios
\be
\frac{1-F_n(x)}{1-\Phi(x)}\quad\text{when}\quad x\ra +\infty,\tag{4a}
\label{eqen4a}
\ee 
and
\be
\frac{F_n(x)}{\Phi(x)}\quad\text{when}\quad x\ra -\infty.\tag{4b}
\label{eqen4b}
\ee
If $x$ is independent of $n$, it follows from \eqref{eqen3} that these ratios both converge to 1 when $n$ goes to infinity; it now remains to understand what happens when $x$ becomes infinite with $n$.

\switchcolumn*

On sait que le théorème de Liapounoff fournit, sous certaines conditions, une borne supérieure du module $| F_n(x)-\Phi(x)|$ qui est du même ordre de grandeur que $\frac{\log n}{\sqrt{n}}$ (voir le chap.\ III). Quand ce théorème est applicable, on trouve sans difficulté que les rapports (\hyperref[eqfr4a]{4}) tendent encore vers l'unité lorsque $|x|$ reste inférieur à $(\frac{1}{2}-\epsilon)\sqrt{\log n}$, où $\epsilon >0$. Cependant ce résultat semble bien insuffisant, notre problème étant d'étudier le comportement des rapports (\hyperref[eqfr4a]{4}) dans un domaine beaucoup plus étendu, par exemple pour des valeurs de $|x|$ qui sont du même ordre de grandeur qu'une puissance de $n$.

\switchcolumn

We know that the theorem of Lyapunov provides, under certain conditions, an upper bound on the modulus $|F_n(x)-\Phi(x)|$ which is of the same order of magnitude as $\frac{\log n}{\sqrt{n}}$ (see Chap.\ III). When this theorem is applicable, we find without difficulty that the ratios (\hyperref[eqen4a]{4}) still converge to 1 when $|x|$ remains below $(\frac{1}{2}-\epsilon)\sqrt{\log n}$, where $\epsilon>0$. However, this result seems well insufficient, our problem being to study the behaviour of the ratios (\hyperref[eqen4a]{4}) in a much larger domain, for example, for values of $|x|$ that are of the same order of magnitude as a power of $n$. \endnote{This is the main problem to be studied: to treat fluctuations of the sum of random variables that are of order $n$, compared to order $\sqrt{n}$ in the central limit theorem.}

\switchcolumn*

Avant d'aborder ce problème, observons qu'on ne peut guère espérer \emph{a priori} d'obtenir des expressions asymptotiques à la fois simples et générales qu'en se bornant aux valeurs de $x$ qui sont de la forme $o(\sqrt{n})$. En effet, la fonction de répartition $V(x)$, qui représente les données du problème, peut être choisie de manière que toute sa variation reste comprise dans un intervalle fini $(-\mu\sigma, \mu\sigma)$. Pour la fonction $F_n(x)$, toute la variation sera alors comprise dans l'intervalle $(-\mu\sqrt{n},\mu\sqrt{n})$, ce qui montre que les rapports (\hyperref[eqfr4a]{4}) s'annuleront identiquement pour $x>\mu\sqrt{n}$ et pour $x<-\mu\sqrt{n}$ respectivement.

\switchcolumn

Before tackling this problem, observe that we can hardly hope \emph{a priori} to obtain asymptotic expressions that are both simple and general only by considering values of $x$ of the form $o(\sqrt{n})$. Indeed, the cumulative distribution function $V(x)$, which represents the data of the problem, can be chosen such that its variation is contained in a finite interval $(\mu\sigma,\mu\sigma)$. For the function $F_n(x)$, all its variation will thus be contained in the interval $(-\mu\sqrt{n},\mu\sqrt{n})$, which shows that the ratios \mbox{(\hyperref[eqen4a]{4})} will cancel identically for $x>\mu\sqrt{n}$ and for $x<-\mu\sqrt{n}$, respectively.

\switchcolumn*

Dans ce qui va suivre, nous allons étudier le comportement asymptotique des rapports (\hyperref[eqfr4a]{4}) en imposant toujours à la fonction $V(x)$ la \emph{condition~\ref{condfrA}} qui va être formulée à l'instant, et en supposant $x$ de la forme $o(\frac{\sqrt{n}}{\log n})$. On verra cependant plus tard (Chap.\ V) que, si la fonction $V(x)$ satisfait à une certaine condition additionnelle~\ref{condfrB}, on peut même considérer les valeurs de $x$ qui sont du même ordre de grandeur que $\sqrt{n}$.

\switchcolumn

In what follows, we shall study the asymptotic behavior of the ratios (\hyperref[eqen4a]{4}) by always imposing to the function $V(x)$ the \emph{condition}~\ref{condenA} formulated just now, and by assuming $x$ of the form $o(\frac{\sqrt{n}}{\log n})$. However, we will see later (Chap.\ V) that, if the function $V(x)$ satisfies a certain additional condition~\ref{condenB}, we can even consider values of $x$ that are of the same order of magnitude as $\sqrt{n}$.~\endnote{This is the large deviation regime of fluctuations.}

\switchcolumn*

\setcounter{equation}{4}
\begin{conditionth}
\label{condfrA}
Il existe un nombre $A>0$ tel que l'intégrale
\be
R=\int_{-\infty}^\infty e^{hy}\, dV(y)
\label{eqfr5}
\ee
converge pour $|h|<A$.
\end{conditionth}

\switchcolumn

\setcounter{equation}{4}
\begin{conditionth}
\label{condenA}
There exists a number $A>0$ such that the integral
\be
R=\int_{-\infty}^\infty e^{hy}\, dV(y)
\label{eqen5}
\ee
converges for $|h|<A$. \emph{\endnote{Condition~\ref{condenA} is now called the Cram\'er condition; it requires the generating function of the CDF $V(x)$ to exist in a neighbourhood of $h=0$.}}
\end{conditionth}

\switchcolumn*

En supposant que cette condition soit satisfaite, nous allons établir entre autres le résultat fondamental suivant. \emph{Pour $x>1$, $x=o(\frac{\sqrt{n}}{\log n})$, on a}
\[
\frac{1-F_n(x)}{1-\Phi(x)}=e^{\frac{x^3}{\sqrt{n}} \lambda\left(\frac{x}{\sqrt{n}}\right)}\left[ 1+O\left(\frac{x \log n}{\sqrt{n}}\right)\right],
\]
\[
\frac{F_n(-x)}{\Phi(-x)}=e^{-\frac{x^3}{\sqrt{n}} \lambda\left(-\frac{x}{\sqrt{n}}\right)}\left[ 1+O\left(\frac{x \log n}{\sqrt{n}}\right)\right],
\]
où
\[
\lambda(z)=\sum_{0}^\infty c_\nu z^\nu
\]
\emph{est une série de puissances convergente pour toute valeur suffisamment petite de $|z|$.}

\switchcolumn

Assuming that this condition is satisfied, we shall establish the following fundamental results among others. \textit{For $x>1$, $x=o(\frac{\sqrt{n}}{\log n})$, we have}
\[
\frac{1-F_n(x)}{1-\Phi(x)}=e^{\frac{x^3}{\sqrt{n}} \lambda\left(\frac{x}{\sqrt{n}}\right)}\left[ 1+O\left(\frac{x \log n}{\sqrt{n}}\right)\right],
\]
\[
\frac{F_n(-x)}{\Phi(-x)}=e^{-\frac{x^3}{\sqrt{n}} \lambda\left(-\frac{x}{\sqrt{n}}\right)}\left[ 1+O\left(\frac{x \log n}{\sqrt{n}}\right)\right],
\]
where
\[
\lambda(z)=\sum_{0}^\infty c_\nu z^\nu
\]
\emph{is a power series convergent for small enough values of $|z|$.}

\switchcolumn*

Ce théorème, dont nous déduirons plusieurs corollaires importants, sera démontré dans le chapitre III. Par l'introduction de la condition additionnelle~\ref{condfrB}, nous parviendrons dans le chapitre V à des résultats encore plus précis. Enfin, dans le dernier chapitre, nous donnerons des théorèmes analogues aux précédents pour le cas d'un processus stocastique homogène.

\switchcolumn

This theorem, from which we will deduce many important corollaries, will be proved in Chap.\ III. By introducing the additional condition~\ref{condenB}, we will reach in Chap.\ V even more precise results. Finally, in the last chapter, we shall give analogous theorems for the case of a homogeneous stochastic process.~\endnote{The original article contains ``stocastique'', which seems to be an error deriving from the Italian ``stocastico''. }

\switchcolumn*

\section*{Chapitre II}

\switchcolumn

\section*{Chapter II}

\switchcolumn*

Supposons que la condition~\ref{condfrA} soit satisfaite et choisissons un nombre réel $h$ situé dans le domaine de convergence de l'intégrale \eqref{eqfr5}. On s’assure facilement que la fonction (\footnote{Une transformation analogue a été employée par F.\ Esscher\,: On the probability function in the collective theory of risk, \emph{Skandinavisk Akzuarietidskrift}, 15 (1932), p.\ 175.})
\be
\bV(x)=\frac{1}{R}\int_{-\infty}^x e^{hy} \, dV(y)
\label{eqfr6}
\ee
possède toutes les propriétés essentielles d’une fonction de répartition.

\switchcolumn

Suppose that Condition~\ref{condenA} is satisfied and choose a real number $h$ in the domain of convergence of the integral \eqref{eqen5}. We can easily verify that the function (\footnote{An analogous transformation has been used by F. Esscher: On the probability function in the collective theory of risk, \textit{Skandinavisk Aktuarietidskrift}, 15 (1932), p.\ 175.}) \endnote{Footnotes are numbered consecutively rather than by page, as in the French original.}
\be
\bV(x)=\frac{1}{R}\int_{-\infty}^x e^{hy} \, dV(y)
\label{eqen6}
\ee
possesses all the essential properties of a cumulative distribution function.~\endnote{The transformation from $V(x)$ to $\bV(x)$, introduced by Esscher in 1932, is now often called the Cram\'er transform.}

\switchcolumn*

Considérons donc une suite $\bZ_1,\bZ_2,\ldots$ de variables aléa\-toires indépendantes ayant toutes la même fonction de réparti\-tion $\bV(x)$. Posons
\be
E(\bZ_n)=\bm,\qquad E(\bZ^2_n)=\bsigma^2.
\label{eqfr7}
\ee

\switchcolumn

Consider then a sequence $\bZ_1,\bZ_2,\ldots$ of independent random variables all having the same cumulative distribution function $\bV(x)$. Let \endnote{The bars over various symbols are often positioned incorrectly in the French original; the Russian version is slightly better.}
\be
E(\bZ_n)=\bm,\qquad E(\bZ^2_n)=\bsigma^2.
\label{eqen7}
\ee

\switchcolumn*

Désignons encore par $\bW(x)$ la fonction de répartition de la somme
\[
\bZ_1+\cdots+\bZ_n,
\]
et par $\bF_n(x)$ la fonction de répartition de la variable
\[
\frac{\bZ_1+\cdots +\bZ_n-\bm n}{\bsigma\sqrt{n}}.
\]

\switchcolumn

Denote further by $\bW_n(x)$ the cumulative distribution function of the sum
\[
\bZ_1+\cdots+\bZ_n,
\]
and by $\bF_n(x)$ the cumulative distribution function of the variable
\[
\frac{\bZ_1+\cdots +\bZ_n-\bm n}{\bsigma\sqrt{n}}.
\]

\switchcolumn*

On a alors
\be
\bF_n(x)=\bW_n(\bm n+\bsigma x\sqrt{n}).
\label{eqfr8}
\ee

\switchcolumn

Hence we have
\be
\bF_n(x)=\bW_n(\bm n+\bsigma x\sqrt{n}).
\label{eqen8}
\ee

\switchcolumn*

Le but de ce chapitre est d’établir une relation entre les fonctions $F_n(x)$ and $\bF_n(x)$ définies par \eqref{eqfr2} et \eqref{eqfr8} respectivement. 

\switchcolumn

The goal of this chapter is to establish a relation between the functions $F_n(x)$ and $\bF_n(x)$ defined by \eqref{eqen2} and \eqref{eqen8}, respectively.

\switchcolumn*

Introduisons les fonctions caractéristiques $\nu(z)$, $\bnu(z)$, $w_n(z)$ et $\bw_n(z)$, où
\[
\nu(z)=\int_{-\infty}^\infty e^{izy}\, dV(y),
\]
tandis que $\bnu(z)$, $w_n(z)$ et $\bw_n(z)$ sont définies par des expressions analogues où figurent les fonctions $\bV$, $W_n$ et $\bW_n$, respectivement. Considérons ici $z$ comme une variable complexe et posons $z=\xi+i\eta$. D’après la condition~\ref{condfrA} et la relation \eqref{eqfr6}, la fonction $\nu(z)$ est holomorphe à l’intérieur de la bande $|\eta|<A$, et l'on a
\[
\begin{split}
R 		&=\nu(-ih),\\
\bnu(z)	&= \frac{1}{R}\nu(z-ih).
\end{split}
\]

\switchcolumn

Let us introduce the characteristic functions $\nu(z)$, $\bnu(z)$, $w_n(z)$ and $\bw_n(z)$, where
\[
\nu(z)=\int_{-\infty}^\infty e^{izy}\, dV(y),
\]
while $\bnu(z)$, $w_n(z)$ and $\bw_n(z)$ are defined by analogous expressions involving the functions $\bV$, $W_n$ and $\bW_n$, respectively. Consider $z$ here as a complex variable and let $z=\xi+i\eta$. Following Condition~\ref{condenA} and the relation \eqref{eqen6}, the function $\nu(z)$ is holomorphic inside the strip $|\eta|<A$, and we have
\[
\begin{split}
R&=\nu(-ih),\\
\bnu(z)&= \frac{1}{R}\nu(z-ih).
\end{split}
\]

\switchcolumn*

D’autre part, les variables $Z_n$ ainsi que les variables $\bZ_n$ étant mutuellement indépendantes, il est bien connu qu’on a
\[
\begin{split}
w_n(z)	&= [\nu(z)]^n,\\
\bw_n(z)	&= [\bnu(z)]^n,
\end{split}
\]
d'où
\be
\bw_n(z)=\frac{1}{R^n}w_n(z-ih),
\label{eqfr9}
\ee
les deux membres de cette dernière relation étant des fonctions holomorphes de $z$ à l’intérieur de la bande $|\eta-h|<A$. En posant ici $z=ih$, il vient
\[
\frac{1}{R^n}=\bw_n(ih)=\int_{-\infty}^\infty e^{-hy}\, d\bW_n(y).
\]

\switchcolumn

On the other hand, the variables $Z_n$ as well as the variables $\bZ_n$ being mutually independent, it is well known that we have
\[
\begin{split}
w_n(z)&=[\nu(z)]^n,\\
\bw_n(z)&= [\bnu(z)]^n,
\end{split}
\]
whence
\be
\bw_n(z)=\frac{1}{R^n}w_n(z-ih),
\label{eqen9}
\ee
the two sides of this last relation being holomorphic functions of $z$ inside the strip $|\eta-h|<A$. Letting here $z=ih$, it follows
\[
\frac{1}{R^n}=\bw_n(ih)=\int_{-\infty}^\infty e^{-hy}\, d\bW_n(y).
\]

\switchcolumn*

Remplaçons dans \eqref{eqfr9} $z$ par $z+ ih$\,; nous avons alors
\be
w_n(z)=R^n \bw_n(z+ih)=\frac{\bw_n(z+ih)}{\bw_n(ih)}.
\label{eqfr10}
\ee

\switchcolumn

Let us replace in \eqref{eqen9} $z$ by $z+ih$, we then have
\be
w_n(z)=R^n \bw_n(z+ih)=\frac{\bw_n(z+ih)}{\bw_n(ih)}.
\label{eqen10}
\ee

\switchcolumn*

Or cette relation équivaut à 
\be
W_n(x)=R^n\int_{-\infty}^x e^{-hy}\, d\bW_n(y).
\label{eqfr11}
\ee
En effet, les deux membres de \eqref{eqfr11} sont des fonctions de répartition en $x$, dont les fonctions caractéristiques, comme on le voit sans difficulté, sont égales à $w_n(z)$ et $R^n\bw_n(z+ih)$ respectivement. Donc, d’après \eqref{eqfr10}, ces fonctions de répartition sont identiques.

\switchcolumn

This relation is equivalent to
\be
W_n(x)=R^n\int_{-\infty}^x e^{-hy}\, d\bW_n(y).
\label{eqen11}
\ee
Indeed, the two members of \eqref{eqen11} are cumulative distribution function in $x$, whose characteristic functions, as we can easily see, are equal to $w_n(z)$ and $R^n\bw_n(z+ih)$, respectively. Therefore from \eqref{eqen10}, these cumulative distribution functions are identical.

\switchcolumn*

En introduisant dans \eqref{eqfr11} les fonctions $F_n(x)$ et $\bF_n(x)$ définies par \eqref{eqfr2} et \eqref{eqfr8}, on obtient au moyen d’une substitution simple
\be
F_n(x)=R^n e^{-h\bm n}\int_{-\infty}^{\frac{\sigma x-\bm\sqrt{n}}{\bsigma}} e^{-h\bsigma\sqrt{n} y}\, d\bF_n(y).\tag{12a}
\label{eqfr12a}
\ee
Si on fait ici tendre $x$ vers l'infini positif on obtient de plus, en faisant la différence,
\be
1-F_n(x)=R^n e^{-h\bm n}\int_{\frac{\sigma x-\bm\sqrt{n}}{\bsigma}}^\infty e^{-h\bsigma\sqrt{n} y}\, d\bF_n(y).\tag{12b}
\label{eqfr12b}
\ee
Les equations \eqref{eqfr12a} et \eqref{eqfr12b} expriment la relation entre $F_n(x)$ et $\bF_n(x)$ que nous nous sommes proposés d’établir.

\switchcolumn

Introducing in \eqref{eqen11} the functions $F_n(x)$ and $\bF_n(x)$ defined by \eqref{eqen2} and \eqref{eqen8}, we obtain by means of a simple substitution
\be
F_n(x)=R^n e^{-h\bm n}\int_{-\infty}^{\frac{\sigma x-\bm\sqrt{n}}{\bsigma}} e^{-h\bsigma\sqrt{n} y}\, d\bF_n(y).\tag{12a}
\label{eqen12a}
\ee
If we let here $x$ converge to positive infinity, we also obtain, by taking the difference,
\be
1-F_n(x)=R^n e^{-h\bm n}\int_{\frac{\sigma x-\bm\sqrt{n}}{\bsigma}}^\infty e^{-h\bsigma\sqrt{n} y}\, d\bF_n(y).\tag{12b}
\label{eqen12b}
\ee
The equations \eqref{eqen12a} and \eqref{eqen12b} express the relation between $F_n(x)$ and $\bF_n(x)$, which we proposed to study. 

\switchcolumn*

Or, d’après la définition \eqref{eqfr8} de $\bF_n(x)$, le théorème de Laplace-Liapounoff peut encore être appliqué à cette fonction. Ce théorème nous apprend que $\bF_n(x)$ tend, pour $n\ra\infty$, vers la fonction de répartition normale $\Phi(x)$ définie par \eqref{eqfr3}. En remplaçant dans \eqref{eqfr12a} et \eqref{eqfr12b} $\bF_n(y)$ par $\Phi(y)$, on doit donc obtenir des expressions approchées de la fonction $F_n(x)$. C'est en précisant ce raisonnement que nous parviendrons, dans les chapitres suivants, à des résultats nouveaux concernant l'allure asymptotique de $F_n(x)$.

\switchcolumn

However, following the definition \eqref{eqen8} of $\bF_n(x)$, the theorem of Laplace-Lyapunov can again be applied to this function. This theorem tells us that $\bF_n(x)$ converges [tends], for $n\ra\infty$, to the normal cumulative distribution function $\Phi(x)$ defined by \eqref{eqen3}. Replacing in \eqref{eqen12a} and \eqref{eqen12b} $\bF_n(y)$ by $\Phi(y)$, we must then obtain expressions close to the function $F_n(x)$. It is by making this reasoning more precise that we will reach, in the following chapters, new results concerning the asymptotic form of $F_n(x)$.

\switchcolumn*

\section*{Chapitre III}

\switchcolumn

\section*{Chapter III}

\switchcolumn*

\setcounter{equation}{12}
Considérons maintenant $h$ comme une variable réelle, et posons
\be
\log R=\log\int_{-\infty}^\infty e^{hy}dV(y)=\sum_{2}^\infty\frac{\gamma_\nu}{\nu!}h^\nu,
\label{eqfr13}
\ee
où, d’après la condition~\ref{condfrA}, la série converge certainement, pour toute valeur suffisamment petite de $|h|$. Les coefficients $\gamma_\nu$ sont les semi-invariants de la fonction de répartition $V(x)$, et il suit de \eqref{eqfr1} qu’on a $\gamma_1=0$, $\gamma_2=\sigma^2$. Par \eqref{eqfr6} et \eqref{eqfr7}, on a de plus
\be
\bm=\frac{d}{dh} \log R=\sum_{2}^\infty \frac{\gamma_\nu}{(\nu-1)!} h^{\nu-1},
\label{eqfr14}
\ee
\be
\bsigma^2=\frac{d\bm}{dh}=\sum_{2}^\infty\frac{\gamma_\nu}{(\nu-1)!}h^{\nu-2}.
\label{eqfr15}
\ee

\switchcolumn

\setcounter{equation}{12}
Consider now $h$ to be a real variable and let
\be
\log R=\log\int_{-\infty}^\infty e^{hy}dV(y)=\sum_{2}^\infty\frac{\gamma_\nu}{\nu!}h^\nu,
\label{eqen13}
\ee
where, following Condition~\ref{condenA}, the series converges certainly for all small enough value $|h|$. The coefficients $\gamma_\nu$ are the semi-invariants \endnote{Cumulants.} of the cumulative distribution function $V(x)$ and it follows from \eqref{eqen1} that we have $\gamma_1=0$, $\gamma_2=\sigma^2$. From \eqref{eqen6} and \eqref{eqen7}, we also have
\be
\bm=\frac{d}{dh} \log R=\sum_{2}^\infty \frac{\gamma_\nu}{(\nu-1)!} h^{\nu-1},
\label{eqen14}
\ee
\be
\bsigma^2=\frac{d\bm}{dh}=\sum_{2}^\infty\frac{\gamma_\nu}{(\nu-1)!}h^{\nu-2}.
\label{eqen15}
\ee

\switchcolumn*

Dans le voisinage du point $h=0$, $\bm$ est donc une fonction continue et toujours croissante de la variable réelle $h$. 

\switchcolumn

In the neighborhood of $h=0$, $\bm$ is thus a continuous and always increasing function of the real variable $h$. 

\switchcolumn*

Soit maintenant $z$ un nombre réel donné. L’équation
\be
\sigma z=\bm=\sum_2^\infty \frac{\gamma_\nu}{(\nu-1)!} h^{\nu-1}
\label{eqfr16}
\ee
admet alors, pour tout $z$ de module suffisamment petit, une seule racine réelle $h$, qui a le même signe que $z$ et tend vers zéro avec $z$. Réciproquement, cette racine $h$ peut être développée en série de puissances de $z$, convergente pour tout $z$ de module suffisamment petit. Les premiers termes de ce développement sont 
\be
h=\frac{z}{\sigma}-\frac{\gamma_3}{2\sigma^4}z^2-\frac{\sigma^2\gamma_4-3\gamma_3^2}{6\sigma^7}z^3+\cdots.
\label{eqfr17}
\ee

\switchcolumn

Now, let $z$ a real number. The equation
\be
\sigma z=\bm=\sum_2^\infty \frac{\gamma_\nu}{(\nu-1)!} h^{\nu-1}
\label{eqen16}
\ee
then admits, for all $z$ of modulus sufficiently small, a unique root $h$, which has the same sign as $z$ and converges to zero with $z$. Conversely, this root $h$ can be expanded in a power series in $z$, which converges for all $x$ of modulus sufficiently small. The first terms of this series are
\be
h=\frac{z}{\sigma}-\frac{\gamma_3}{2\sigma^4}z^2-\frac{\sigma^2\gamma_4-3\gamma_3^2}{6\sigma^7}z^3+\cdots.
\label{eqen17}
\ee

\switchcolumn*

De \eqref{eqfr13} et \eqref{eqfr14} on tire
\[
h\bm-\log R=\sum_2^\infty \frac{(\nu-1)\gamma_\nu}{\nu!} h^\nu.
\]

\switchcolumn

From \eqref{eqen13} and \eqref{eqen14} we get
\[
h\bm-\log R=\sum_2^\infty \frac{(\nu-1)\gamma_\nu}{\nu!} h^\nu.
\]

\switchcolumn*

En remplaçant ici $h$ par son développement \eqref{eqfr17}, il vient
\[
h\bm-\log R=\frac{1}{2}z^2-\frac{\gamma_3}{6\sigma^3}z^3-\frac{\sigma^2\gamma_4-3\gamma_3^2}{24\sigma^6}z^4+\cdots.
\]

\switchcolumn

Replacing here $h$ by its expansion \eqref{eqen17}, it follows
\[
h\bm-\log R=\frac{1}{2}z^2-\frac{\gamma_3}{6\sigma^3}z^3-\frac{\sigma^2\gamma_4-3\gamma_3^2}{24\sigma^6}z^4+\cdots.
\]

\switchcolumn*

Donc, en observant que $z^2=(\frac{\bm}{\sigma})^2$ et en posant 
\be
\frac{\bm^2}{2\sigma^2}-h\bm+\log R=z^3\lambda(z),
\label{eqfr18}
\ee
$\lambda(z)$ admet un développement en série de puissances
\be
\lambda(z)=c_0+c_1 z+c_2 z^2+\cdots
\label{eqfr19}
\ee
convergent pour tout $z$ de module suffisamment petit, et l'on a 
\be
c_0=\frac{\gamma_3}{6\sigma^3},\qquad c_1=\frac{\sigma^2\gamma_4-3\gamma_3^2}{24\sigma^6},\cdots.
\label{eqfr20}
\ee

\switchcolumn

Hence, observing that $z^2=(\frac{\bm}{\sigma})^2$ and letting
\be
\frac{\bm^2}{2\sigma^2}-h\bm+\log R=z^3\lambda(z),
\label{eqen18}
\ee
$\lambda(z)$ admits a series expansion
\be
\lambda(z)=c_0+c_1 z+c_2 z^2+\cdots
\label{eqen19}
\ee
which converges for all $z$ of modulus sufficiently small, and we have 
\be
c_0=\frac{\gamma_3}{6\sigma^3},\qquad c_1=\frac{\sigma^2\gamma_4-3\gamma_3^2}{24\sigma^6},\cdots.
\label{eqen20}
\ee

\switchcolumn*

Ceci posé, nous pouvons énoncer notre théorème fondamental de la manière suivante.

\switchcolumn

This being given, we can formulate our fundamental theorem in the following way.

\switchcolumn*

\begin{theoremfr}
\label{thmfr1}
Supposons que la condition~\ref{condfrA} soit satisfaite. Soit $x$ un nombre réel qui peut dépendre de $n$, tel que $x>1$ et $x=o(\frac{\sqrt{n}}{\log n})$
lorsque $n$ tend vers l'infini. Pour la fonction de répartition $F_n(x)$ introduite dans le chapitre I, on a alors
\[
\frac{1-F_n(x)}{1-\Phi(x)}=e^{\frac{x^3}{\sqrt{n}} \lambda(\frac{x}{\sqrt{n}})}\left[1+O\left(\frac{x\log n}{\sqrt{n}}\right)\right]
\]
et
\[
\frac{F_n(-x)}{\Phi(-x)}=e^{-\frac{x^3}{\sqrt{n}} \lambda(-\frac{x}{\sqrt{n}})}\left[1+O\left(\frac{x\log n}{\sqrt{n}}\right)\right]
\]
$\lambda(z)$ étant la fonction définie par \eqref{eqfr18} et qui admet le dévelop\-pe\-ment \eqref{eqfr19}, dont les premiers coefficients sont donnés par \eqref{eqfr20}.
\end{theoremfr}

\switchcolumn

\begin{theoremen}
\label{thmen1}
Assume that Condition~\ref{condenA} is satisfied. Let $x$ be a real number that can depend on $n$ such that $x>1$ and $x=o(\frac{\sqrt{n}}{\log n})$
when $n$ goes to infinity. For the cumulative distribution function $F_n(x)$ introduced in chapter 1, we then have
\[
\frac{1-F_n(x)}{1-\Phi(x)}=e^{\frac{x^3}{\sqrt{n}} \lambda(\frac{x}{\sqrt{n}})}\left[1+O\left(\frac{x\log n}{\sqrt{n}}\right)\right]
\]
and
\[
\frac{F_n(-x)}{\Phi(-x)}=e^{-\frac{x^3}{\sqrt{n}} \lambda(-\frac{x}{\sqrt{n}})}\left[1+O\left(\frac{x\log n}{\sqrt{n}}\right)\right]
\]
$\lambda(z)$ being the function defined by \eqref{eqen18} and which admits the development \eqref{eqen19}, whose first coefficients are given by \eqref{eqen20}.
\end{theoremen}

\switchcolumn*

Les démonstrations des deux relations énoncées étant tout à fait analogues, nous nous bornerons à la démonstration de la première relation.

\switchcolumn

The proofs of the two relations just stated being totally analogous, we will restrict ourselves to prove the first relation.

\switchcolumn*

Si dans \eqref{eqfr12b} nous prenons $x=\frac{\bm\sqrt{n}}{\sigma}$, nous aurons
\be
1-F_n(\frac{\bm\sqrt{n}}{\sigma})=R^n e^{-h\bm n} \int_0^\infty e^{-h\bsigma\sqrt{n} y} d\bF_n(y)
\label{eqfr21}
\ee
pour toute valeur réelle de $h$ telle que $|h| <A$. Posons maintenant
\be
\bF_n(y)=\Phi(y)+Q_n(y).
\label{eqfr22}
\ee

\switchcolumn

If in \eqref{eqen12b} we take $x=\frac{\bm\sqrt{n}}{\sigma}$, we will have
\be
1-F_n(\frac{\bm\sqrt{n}}{\sigma})=R^n e^{-h\bm n} \int_0^\infty e^{-h\bsigma\sqrt{n} y} d\bF_n(y)
\label{eqen21}
\ee
for all real values of $h$ such that $|h| <A$. Now let
\be
\bF_n(y)=\Phi(y)+Q_n(y).
\label{eqen22}
\ee

\switchcolumn*

D’après le théorème de Liapounoff, on a (\footnote{Voir H. Cramér, Random variables and probability distributions, \emph{Cambridge Tracts in Mathematics}, No 36, Cambridge 1937, p.\ 77.}) alors pour tout $n > 1$ et pour tout $y$ réel
\[
|Q_n(y)|<k\frac{\log n}{\sqrt{n}},
\]
où
\[
k=\frac{3}{\bsigma^3}\int_{-\infty}^\infty |y-\bm|^3 dV(y).
\]

\switchcolumn

Following the theorem of Liapounov, we then have (\footnote{See H. Cramér, Random variables and probability distributions, Cambridge Tracts in Mathematics, No 36, Cambridge 1937, p.\ 77.}) for all $n > 1$ and for all real $y$
\[
|Q_n(y)|<k\frac{\log n}{\sqrt{n}},
\]
where
\[
k=\frac{3}{\bsigma^3}\int_{-\infty}^\infty |y-\bm|^3 dV(y).
\]

\switchcolumn*

Le nombre $k$ ainsi défini dépend évidemment de $h$. Or il suit de \eqref{eqfr14}, \eqref{eqfr15} et de la condition~\ref{condfrA} qu’on peut déterminer un nombre positif $A_1< A$ tel que pour $|h|<A_1$, on ait $k<K$ et par conséquent
\be
|Q_n(y)|<K \frac{\log n}{\sqrt{n}},
\label{eqfr23}
\ee
où $K$ est indépendant de $h$, $n$ et $y$. 

\switchcolumn

The number $k$ thus defined depends obviously on $h$. However, it follows from \eqref{eqen14}, \eqref{eqen15} and Condition~\ref{condenA} that we can determine a positive number $A_1< A$ such that for $|h|<A_1$, we have $k<K$ and, consequently,
\be
|Q_n(y)|<K \frac{\log n}{\sqrt{n}},
\label{eqen23}
\ee
where $K$ is independent of $h$, $n$ and $y$.

\switchcolumn*

Dès maintenant, nous considérons $h$ comme une variable essentiellement \textit{positive}. Faisons tendre $n$ vers l'infini et $h$ vers zéro, de manière que le produit ait une borne inférieure positive.

\switchcolumn

From now on, we shall consider $h$ as an essentially \emph{positive} variable. Let $n$ go to infinity and $h$ to zero, in such a way that the product has a positive lower bound.

\switchcolumn*

Nous avons alors par \eqref{eqfr14} et \eqref{eqfr15}
\[
h\bsigma\sqrt{n}=\frac{\bm\sqrt{n}}{\sigma}+O(h^2\sqrt{n}).
\]

\switchcolumn

We then have from \eqref{eqen14} and \eqref{eqen15}
\[
h\bsigma\sqrt{n}=\frac{\bm\sqrt{n}}{\sigma}+O(h^2\sqrt{n}).
\]

\switchcolumn*

De \eqref{eqfr22} et \eqref{eqfr23} on déduit au moyen de calculs faciles
\[
\begin{split}
\int_0^\infty & e^{-h\bsigma\sqrt{n} y}d\bF_n(y) \nonumber\\
&= \frac{1}{2\pi}\int_0^\infty e^{-h\bsigma\sqrt{n} y-\frac{1}{2}y^2}dy-Q_n(0)\\
& \qquad +h\bsigma\sqrt{n}\int_0^\infty e^{-h\bsigma\sqrt{n} y}Q_n(y)dy\nonumber\\
&= \frac{1}{2\pi}\int_0^\infty e^{-h\bsigma\sqrt{n} y-\frac{1}{2}y^2}dy+O\left(\frac{\log n}{\sqrt{n}}\right)\nonumber\\
&= \frac{1}{2\pi}\int_0^\infty e^{-h\bsigma\sqrt{n} y-\frac{1}{2}y^2}dy\cdot [1+O(h\log n)]\nonumber\\
&= \frac{1}{2\pi}\int_0^\infty e^{-\frac{\bm\sqrt{n}}{\sigma} y-\frac{1}{2}y^2}dy\cdot [1+O(h\log n)]\nonumber\\
&= e^{\frac{n{\bm}^2}{2\sigma^2}}\left[1-\Phi\left(\frac{\bm\sqrt{n}}{\sigma}\right)\right] \cdot [1+O(h\log n)].
\end{split}
\]

\switchcolumn

From \eqref{eqen22} and \eqref{eqen23} we deduce through simple calculations
\[
\begin{split}
\int_0^\infty & e^{-h\bsigma\sqrt{n} y}d\bF_n(y) \nonumber\\
&= \frac{1}{2\pi}\int_0^\infty e^{-h\bsigma\sqrt{n} y-\frac{1}{2}y^2}dy-Q_n(0)\\
& \qquad +h\bsigma\sqrt{n}\int_0^\infty e^{-h\bsigma\sqrt{n} y}Q_n(y)dy\nonumber\\
&= \frac{1}{2\pi}\int_0^\infty e^{-h\bsigma\sqrt{n} y-\frac{1}{2}y^2}dy+O\left(\frac{\log n}{\sqrt{n}}\right)\nonumber\\
&= \frac{1}{2\pi}\int_0^\infty e^{-h\bsigma\sqrt{n} y-\frac{1}{2}y^2}dy\cdot [1+O(h\log n)]\nonumber\\
&= \frac{1}{2\pi}\int_0^\infty e^{-\frac{\bm\sqrt{n}}{\sigma} y-\frac{1}{2}y^2}dy\cdot [1+O(h\log n)]\nonumber\\
&= e^{\frac{n{\bm}^2}{2\sigma^2}}\left[1-\Phi\left(\frac{\bm\sqrt{n}}{\sigma}\right)\right] \cdot [1+O(h\log n)].
\end{split}
\]

\switchcolumn*

En introduisant la dernière expression dans \eqref{eqfr21}, on aura
\be
\frac{1-F_n(\frac{\bm\sqrt{n}}{\sigma})}{1-\Phi(\frac{\bm\sqrt{n}}{\sigma})}=e^{n(\frac{\bm^2}{2\sigma^2}-h\bm+\log R)}[1+O(h\log n)].
\label{eqfr24}
\ee

\switchcolumn

Introducing the last expression in \eqref{eqen21}, we will have
\be
\frac{1-F_n(\frac{\bm\sqrt{n}}{\sigma})}{1-\Phi(\frac{\bm\sqrt{n}}{\sigma})}=e^{n(\frac{\bm^2}{2\sigma^2}-h\bm+\log R)}[1+O(h\log n)].
\label{eqen24}
\ee

\switchcolumn*

Soit maintenant $x$ un nombre réel qui peut dépendre de $n$, tel que $x>1$ et $x=o(\frac{\sqrt{n}}{\log n})$. Formons l'équation
\be
x=\frac{\bm\sqrt{n}}{\sigma},
\label{eqfr25}
\ee
qui peut aussi s'écrire
\[
\frac{\sigma x}{\sqrt{n}}=\bm=\sum_{2}^\infty\frac{\gamma_\nu}{(\nu-1)!} h^{\nu-1},
\]
et qui, par la substitution $z=\frac{x}{\sqrt{n}}$, devient identique à l'équation \eqref{eqfr16}. Pour tout $n$ suffisamment grand, l’équation \eqref{eqfr25} admet donc une seule racine positive $h$ qui tend vers zéro lorsque $n$ tend vers l’infini. D’après \eqref{eqfr18} on a
\[
\frac{\bm^2}{2\sigma^2}-h\bm+\log R=\left(\frac{x}{\sqrt{n}}\right)^3\lambda\left(\frac{x}{\sqrt{n}}\right),
\]
où $\lambda(z)$ est défini par \eqref{eqfr19}-\eqref{eqfr20}. Le produit $h\sqrt{n}$ est bien borné inférieurement, car on déduit de \eqref{eqfr17}
\[
h\sim \frac{z}{\sigma}=\frac{x}{\sigma\sqrt{n}},
\]
et nous avons supposé $x > 1$.

\switchcolumn

Now let $x$ be a real number which can depend on $n$, such that $x>1$ and $x=o(\frac{\sqrt{n}}{\log n})$. Let us form the equation
\be
x=\frac{\bm\sqrt{n}}{\sigma},
\label{eqen25}
\ee
which can also be written as
\[
\frac{\sigma x}{\sqrt{n}}=\bm=\sum_{2}^\infty\frac{\gamma_\nu}{(\nu-1)!} h^{\nu-1},
\]
and which, with the substitution $z=\frac{x}{\sqrt{n}}$, becomes identical to equation \eqref{eqen16}. For all $n$ sufficiently large, equation \eqref{eqen25} then admits a unique positive root $h$ which converges to zero when $n$ goes to infinity. From \eqref{eqen18} we have
\[
\frac{\bm^2}{2\sigma^2}-h\bm+\log R=\left(\frac{x}{\sqrt{n}}\right)^3\lambda\left(\frac{x}{\sqrt{n}}\right),
\]
where $\lambda(z)$ is defined by \eqref{eqen19}-\eqref{eqen20}. The product $h\sqrt{n}$ is indeed bounded below, since we deduce from \eqref{eqen17}
\[
h\sim \frac{z}{\sigma}=\frac{x}{\sigma\sqrt{n}},
\]
and we have assumed $x > 1$.

\switchcolumn*

Dans \eqref{eqfr24}, on peut donc prendre $h$ égal à la racine de \eqref{eqfr25}. On obtient ainsi
\[
\frac{1-F_n(x)}{1-\Phi(x)}=e^{\frac{x^3}{\sqrt{n}}\lambda(\frac{x}{\sqrt{n}})}\left[1+O\left(\frac{x\log n}{\sqrt{n}}\right)\right],
\]
et le théorème~\ref{thmfr1} est démontré.

\switchcolumn

In \eqref{eqen24}, we can therefore take $h$ equal to the root of \eqref{eqen25}. We then obtain
\[
\frac{1-F_n(x)}{1-\Phi(x)}=e^{\frac{x^3}{\sqrt{n}}\lambda(\frac{x}{\sqrt{n}})}\left[1+O\left(\frac{x\log n}{\sqrt{n}}\right)\right],
\]
and Theorem~\ref{thmen1} is proved.

\switchcolumn*

Du théorème~\ref{thmfr1}, on peut déduire plusieurs corollaires intéressants. Démontrons d’abord le théorème suivant.

\switchcolumn

From Theorem~\ref{thmen1}, we can deduce several interesting corollaries. Let us first prove the following theorem.

\switchcolumn*

\begin{theoremfr}
\label{thmfr2}
Si la condition~\ref{condfrA} est satisfaite, on a pour $x>1$, $x=O(n^{\frac{1}{6}})$,
\[
\frac{1-F_n(x)}{1-\Phi(x)}= e^{\frac{c_0x^3}{\sqrt{n}}}+O\left(\frac{x\log n}{\sqrt{n}}\right),
\]
\[
\frac{F_n(-x)}{\Phi(-x)}=e^{-\frac{c_0x^3}{\sqrt{n}}}+O\left(\frac{x\log n}{\sqrt{n}}\right).
\]
\end{theoremfr}

\switchcolumn

\begin{theoremen}
\label{thmen2}
If Condition~\ref{condenA} is satisfied, we have for $x>1$, $x=O(n^{\frac{1}{6}})$,
\[
\frac{1-F_n(x)}{1-\Phi(x)}= e^{\frac{c_0x^3}{\sqrt{n}}}+O\left(\frac{x\log n}{\sqrt{n}}\right),
\]
\[
\frac{F_n(-x)}{\Phi(-x)}=e^{-\frac{c_0x^3}{\sqrt{n}}}+O\left(\frac{x\log n}{\sqrt{n}}\right).
\]
\end{theoremen}

\switchcolumn*

Ceci est une conséquence immédiate du théorème~\ref{thmfr1}, si l'on remarque que, pour $x=O(n^{1/6})$, l'exposant $\frac{x^3}{\sqrt{n}}\lambda(\frac{x}{\sqrt{n}})$ reste borné lorsque $n$ tend vers l’infini. On voit en particulier que, si $x=o(n^{1/6})$, les deux rapports considérés tendent vers l'unité quand $n$ tend vers l'infini.

\switchcolumn

This is an immediate consequence of Theorem~\ref{thmen1}, if we remark that, for $x=O(n^{1/6})$, the exponent $\frac{x^3}{\sqrt{n}}\lambda(\frac{x}{\sqrt{n}})$ remains bounded when $n$ goes to infinity. In particular, we see that, if $x=o(n^{1/6})$, the two ratios considered converge to 1 when $n$ goes to infinity.

\switchcolumn*

En observant que l'on a pour $x > 1$,
\[
1-\Phi(x)<\frac{1}{x\sqrt{2\pi}}e^{-\frac{x^2}{2}},
\]
\[
\Phi(-x)<\frac{1}{x\sqrt{2\pi}}e^{-\frac{x^2}{2}},
\]
on obtient aussi sans difficulté le théorème suivant qui se rattache immédiatement au théorème de Liapounoff.

\switchcolumn

By observing that we have for $x > 1$,
\[
1-\Phi(x)<\frac{1}{x\sqrt{2\pi}}e^{-\frac{x^2}{2}},
\]
\[
\Phi(-x)<\frac{1}{x\sqrt{2\pi}}e^{-\frac{x^2}{2}},
\]
we also obtain without difficulty the following theorem which immediately relates to the theorem of Liapounov.

\switchcolumn*

\begin{theoremfr}
\label{thmfr3}
Si la condition~\ref{condfrA} est satisfaite, on a pour $x>0$, $x=O(n^{1/6})$,
\[
1-F_n(x)=[1-\Phi(x)]e^{\frac{c_0 x^3}{\sqrt{n}}}+O\left(\frac{\log n}{\sqrt{n}}e^{-\frac{x^2}{2}}\right),
\]
\[
F_n(-x)=\Phi(-x)e^{-\frac{c_0 x^3}{\sqrt{n}}}+O\left(\frac{\log n}{\sqrt{n}}e^{-\frac{x^2}{2}}\right).
\]
\end{theoremfr}

\switchcolumn

\begin{theoremen}
\label{thmen3}
If Condition~\ref{condenA} is satisfied, we have for $x>0$, $x=O(n^{1/6})$,
\[
1-F_n(x)=[1-\Phi(x)]e^{\frac{c_0 x^3}{\sqrt{n}}}+O\left(\frac{\log n}{\sqrt{n}}e^{-\frac{x^2}{2}}\right),
\]
\[
F_n(-x)=\Phi(-x)e^{-\frac{c_0 x^3}{\sqrt{n}}}+O\left(\frac{\log n}{\sqrt{n}}e^{-\frac{x^2}{2}}\right).
\]
\end{theoremen}

\switchcolumn*

Si, tout en restant dans les conditions du théorème~\ref{thmfr1}, $x$ est d’un ordre de grandeur plus élevé que celui de $n^{1/6}$, le théorème~\ref{thmfr1} fournit encore des expressions asymptotiques des probabilités $1-F_n(x)$ et $F_n(-x)$. Si, par exemple, le coefficient $c_0=\frac{\gamma_3}{6\sigma^3}$ est différent de zéro, on voit que l'exposant $\frac{x^3}{\sqrt{n}}\lambda(\frac{x}{\sqrt{n}})$ est asymptotiquement équivalent à $\frac{c_0 x^3}{\sqrt{n}}$. Lorsque $\frac{x}{n^{1/6}}$ tend vers l'infini positif, cet exposant tend donc vers $+\infty$ ou vers $-\infty$ selon le signe de $c_0$. D'après le théorème~\ref{thmfr1}, le rapport $\frac{1-F_n(x)}{1-\Phi(x)}$ tend vers $+\infty$ ou vers $0$ suivant le cas. Par le même raisonnement, le rapport $\frac{F_n(-x)}{\Phi(-x)}$ tend alors vers $0$ ou vers $+\infty$ respectivement.~--- Si $c_0=0$, c’est évidemment le premier coefficient $c_\nu\neq 0$ qui va dominer la question, sans qu’il soit nécessaire d'en préciser ici tous les détails.

\switchcolumn

If, staying in the conditions of Theorem~\ref{thmen1}, $x$ is of order higher than the order of $n^{1/6}$, Theorem~\ref{thmen1} still provides asymptotic expressions for the probabilities $1-F_n(x)$ and $F_n(-x)$. If, for example, the coefficient $c_0=\frac{\gamma_3}{6\sigma^3}$ \endnote{The original article contains $\gamma^3$; this error is repeated in the Russian version.} is different from zero, then we see that the exponent $\frac{x^3}{\sqrt{n}}\lambda(\frac{x}{\sqrt{n}})$ is asymptotically equivalent to $\frac{c_0 x^3}{\sqrt{n}}$. When $\frac{x}{n^{1/6}}$ goes to infinity, this exponent then goes to $+\infty$ or to $-\infty$ depending on the sign of $c_0$. Following Theorem~\ref{thmen1}, the ratio $\frac{1-F_n(x)}{1-\Phi(x)}$ goes to $+\infty$ or to $0$, respectively, depending on the case. By the same reasoning, the ratio $\frac{F_n(-x)}{\Phi(-x)}$ then goes to $0$ or to $+\infty$, respectively.~---~If $c_0=0$, it is obvious that the first coefficient $c_\nu\neq 0$ will dominate the question, so it is not necessary to precise here all the details.

\switchcolumn*

En dernier lieu, un calcul simple permet de déduire du théorème~\ref{thmfr1} Ia généralisation suivante d’un théorème dû à M.\ Khintchine (cf.\ le chapitre suivant).

\switchcolumn

Finally, a simple calculation allows to deduce from Theorem~\ref{thmen1} the following generalization of a theorem due to Khintchine (see next chapter).

\switchcolumn*

\begin{theoremfr}
\label{thmfr4}
Soit $c$ une constante positive. Si la condition~\ref{condfrA} est satisfaite, les deux expressions
\[
\frac{F_n(x+\frac{c}{x})-F_n(x)}{1-F_n(x)}\quad\textrm{et}\quad \frac{\Phi(x+\frac{c}{x})-\Phi(x)}{1-\Phi(x)},
\]
tendent, pour $n\ra\infty$, $x\ra\infty$, $x=O(\frac{\sqrt{n}}{\log n})$, vers une même limite, à savoir vers la quantité $1-e^{-c}$.
\end{theoremfr}

\switchcolumn

\begin{theoremen}
\label{thmen4}
Let $c$ be a positive constant. If Condition~\ref{condenA} is satisfied, the two expressions
\[
\frac{F_n(x+\frac{c}{x})-F_n(x)}{1-F_n(x)}\quad\textrm{and}\quad \frac{\Phi(x+\frac{c}{x})-\Phi(x)}{1-\Phi(x)},
\]
converge, for $n\ra\infty$, $x\ra\infty$, $x=O(\frac{\sqrt{n}}{\log n})$ \emph{\endnote{The original article contains $0$ instead of $O$. The correct $O$ is found in the Russian translation.}}, to the same limit, namely, the quantity $1-e^{-c}$.
\end{theoremen}

\switchcolumn*

Il y a évidemment un théorème correspondant pour les valeurs négatives de la variable.

\switchcolumn

There is evidently a corresponding theorem for the negative values of the variable.

\switchcolumn*

\section*{Chapitre IV}

\switchcolumn

\section*{Chapter IV}

\switchcolumn*

Si, en particulier, on choisit les variables aléatoires $Z_n$ introduites au début de ce travail telles que, pour chaque $Z_n$, il n’y ait que deux valeurs possibles\,:
\[
Z_n=
\left\{
\begin{array}{rll}
1-p & \textrm{avec la probabilité} & p, \\
-p & \quad\textrm{''}\quad\textrm{''}\qquad\textrm{''}& q=1-p,
\end{array}
\right.
\]
on voit qu'on arrive au \textit{cas des épreuves répétées}. On sait que, dans ce cas particulier, la fonction de répartition $F_n(x)$ peut être interprétée de la manière suivante.

\switchcolumn

If, in particular, we choose the random variables $Z_n$ introduced at the beginning of this work such that, for each $Z_n$, there are only two possible values:
\[
Z_n=
\left\{
\begin{array}{rll}
1-p & \textrm{with probability} & p, \\
-p & \quad\textrm{''}\qquad\textrm{''}& q=1-p,
\end{array}
\right.
\]
we see that we arrive at the \emph{case of repeated trials}. We know that, in this particular case, the cumulative distribution function $F_n(x)$ can be interpreted as follows. 

\switchcolumn*

Supposons qu’on fasse une série de $n$ tirages indépendants d'une urne, la probabilité d’amener une boule blanche étant toujours égale à $p$. Désignons par $\nu$ le nombre des boules blanches obtenues au cours des $n$ tirages. Alors nous avons
\[
F_n(x)=\textrm{Prob}(\nu\leq np+x\sqrt{npq}),
\]
et il bien connu que
\[
\lim_{n\ra\infty} F_n(x)=\Phi(x),
\]
pour tout $x$ réel indépendant de $n$.

\switchcolumn

Suppose that we have a series of $n$ independent trials, the probability of choosing a white ball being always equal to $p$. Let $\nu$ be the number of white balls obtained in these $n$ trials. Then we have
\[
F_n(x)=\textrm{Prob}(\nu\leq np+x\sqrt{npq}),
\]
and it is well known that
\[
\lim_{n\ra\infty} F_n(x)=\Phi(x),
\]
for all $x$ real and independent of $n$.

\switchcolumn*

La quantité $\frac{\nu-np}{\sqrt{npq}}$ est (avec un léger changement formel) ce que, d’après M.\ Borel, on appelle l’\textit{écart relatif}. Pour tout $x$ réel indépendant de $n$, la probabilité d’avoir un écart relatif $\leq x$ tend donc vers la limite $\Phi(x)$ lorsque $n$ tend vers l’infini.

\switchcolumn

The quantity $\frac{\nu-np}{\sqrt{npq}}$ is (with a slight formal change) what, following Borel, we call the \textit{relative discrepancy}. For all $x$ real and independent of $n$, the probability of having a relative discrepancy $\leq x$ thus converges to $\Phi(x)$ when $n$ goes to $\infty$.

\switchcolumn*

Cependant, il peut souvent être important de connaître le comportement asymptotique des \textit{probabilités des grands écarts relatifs}, c'est-à-dire le comportement de $F_n(x)$ quand $x$ peut varier avec $n$, en tendant vers $+\infty$, ou vers $-\infty$ quand $n$ croît indéfiniment. Ce cas particulier du problème qui nous occupe dans ce travail a été considéré par plusieurs auteurs (\footnote{Voir p.\ ex.\ N. Smirnoff, Uber Wahrscheinlichkeiten grosser Abweichungen, \textit{Rec.~Soc.~Math.}\ Moscou, 40 (1933), p.\ 441\,; A. Khintchine, Uber einen neuen Grenzwertstaz der Wahrscheinlichkeitsrechnung, \textit{Math.~Annalen} 101 (1929), p.\ 745\,; M.\ Fréchet, Recherches théoriques modernes sur le calcul des probabilités, Paris, 1937, p.\ 222\,; P.\ Lévy, Théorie de l'addition des variables aléatoires, Paris 1937, p.\ 284.\par H. Cramer, etc.}). La plupart des résultats trouvés dans cette direction rentrent dans les théorèmes démontr\'es dans le chapitre précédent.

\switchcolumn

However, it is often important to know the asymptotic behavior of \textit{probabilities of large relative discrepancy}, that is, the behavior of $F_n(x)$ when $x$ can vary with $n$, going to $+\infty$ or to $-\infty$ when $n$ indefinitely grows. This particular case of the problem that we consider in this work has been considered by many authors (\footnote{See e.g.\ N.\ Smirnoff, Uber Wahrscheinlichkeiten grosser Abweichungen, \textit{Rec. Soc. Math.} Moscow, 40 (1933), p.\ 441; A.\ Khintchine, Uber einen neuen Grenzwertstaz der Wahrscheinlichkeitsrechnung, \textit{Math. Annalen} 101 (1929), p.\ 745; M.\ Fréchet, Recherches théoriques modernes sur le calcul des probabilités, Paris, 1937, p.\ 222; P.\ Lévy, Théorie de l'addition des variables aléatoires, Paris 1937, p.\ 284.\par H. Cramer, etc.}). Most results obtained in that direction fit in the theorems proved in the previous chapter.

\switchcolumn*

Ainsi M.\ Smirnoff a démontré un théorème qui peut s’expri\-mer par la relation
\[
\frac{1-F_n(x)}{1-\Phi(x)}=1+o\left(\frac{1}{x^{2s}}\right),1
\]
pour $x=o(n^{\frac{1}{4s+6}})$, $s=0,1,2,\ldots$, et par une relation analogue pour les écarts négatifs. Comme on a pour les valeurs indiquées de $x$
\[
e^{\frac{c_0 x^3}{\sqrt{n}}}+O\left(\frac{x\log n}{\sqrt{n}}\right)=1+O\left(\frac{x^3+x\log n}{\sqrt{n}}\right)=1+o\left(\frac{1}{x^{2s}}\right),
\]
ce résultat est contenu dans notre théorème~\ref{thmfr2}.

\switchcolumn

Thus Smirnoff has proved a theorem which can be expressed by the relation
\[
\frac{1-F_n(x)}{1-\Phi(x)}=1+o\left(\frac{1}{x^{2s}}\right),
\]
for $x=o(n^{\frac{1}{4s+6}})$, $s=0,1,2,\ldots$, and by a similar relation for negative discrepancies. Since we have for the indicated values of $x$
\[
e^{\frac{c_0 x^3}{\sqrt{n}}}+O\left(\frac{x\log n}{\sqrt{n}}\right)=1+O\left(\frac{x^3+x\log n}{\sqrt{n}}\right)=1+o\left(\frac{1}{x^{2s}}\right),
\]
this result is contained in our Theorem~\ref{thmen2}.

\switchcolumn*

D'autre part, M.\ Lévy a donné, pour le cas des épreuves répétées, la relation (sur laquelle nous reviendrons dans le chapitre suivant)
\[
\log(1-F_n(x))\sim \log(1-\Phi(x))\sim -\frac{x^2}{2},
\]
qui est une conséquence de notre théorème~\ref{thmfr1}, et enfin notre théorème~\ref{thmfr4} a été démontré pour le même cas particulier par M.\ Khintchine. Dans ces théorèmes de MM.\ Lévy et Khintchine, notre condition $x=o(\frac{\sqrt{n}}{\log n})$ se trouve remplacée par la condition un peu moins restreinte (\footnote{Dans le cas des épreuves répétées, la fonction de répartition $V(x)$ ne satisfait pas à la condition B, qui va être introduite dans le chapitre suivant et qui nous permettra de remplacer dans nos théorèmes la condition $x=o(\frac{\sqrt{n}}{\log n})$ par la condition $x=o(\sqrt{n})$. Les résultats cités de MM.\ Khintchine et Lévy ne sont donc pas entièrement contenus dans nos résultats.}) $x=o(\sqrt{n})$.

\switchcolumn

Moreover, Lévy has given, for the case of repeated trials, the relation (on which we will come back in the next chapter)
\[
\log(1-F_n(x))\sim \log(1-\Phi(x))\sim -\frac{x^2}{2},
\]
which is a consequence of our Theorem~\ref{thmen1}, and finally our Theorem~\ref{thmen4} has been proved for the same particular case by Khintchine. In those theorems of Lévy and Khintchine, our condition $x=o(\frac{\sqrt{n}}{\log n})$ is replaced by the least restrictive (\footnote{In the case of repeated trials, the cumulative distribution function $V(x)$ does not satisfy Condition~\ref{condenB}, which will be introduced in the next chapter and which will allow us to replace in our theorems the condition $x=o(\frac{\sqrt{n}}{\log n})$ by the condition $x=o(\sqrt{n})$. The cited results of Khintchine and Lévy are thus not entirely contained in our results.}) $x=o(\sqrt{n})$.

\switchcolumn*

\section*{Chapitre V}

\switchcolumn

\section*{Chapter V}

\switchcolumn*

Retournons au problème général posé dans le chapitre I. Jus\-qu'ici, nous avons assujetti la fonction de répartition donnée $V(x)$ à la seule condition~\ref{condfrA}\,; nous allons maintenant introduire une condition additionnelle~\ref{condfrB}, qui nous permettra d'aller plus loin dans l'étude du comportement de la fonction $F_n(x)$ pour des grandes valeurs de $n$ et $x$.

\switchcolumn

Let us go back to the general problem formulated in Chapter~I. Up to now, we have only imposed Condition~\ref{condenA} on the cumulative distribution function $V(x)$; we shall now introduce another condition~\ref{condenB} that will allow us to go further in our study of the behavior of the function $F_n(x)$ for large values of $n$ and $x$.

\switchcolumn*

La fonction $V(x)$ peut, d'une manière unique, être mise sous la forme (\footnote{Voir H. Cramér, \textit{l.\,c.}, p.\ 17.})
\be
V(x)=\beta U_1(x)+(1-\beta)U_2(x)
\label{eqfr26}
\ee
avec $0\leq \beta \leq 1$ où $U_1(x)$ et $U_2(x)$ sont deux fonctions de répartition telles qu'on ait \textit{presque partout}
\[
\begin{split}
U_1(x) &=\int_{-\infty}^x U'_1(y) dy,\\
U'_2(x) &=0.
\end{split}
\]

\switchcolumn

The function $V(x)$ can be put, in a unique way, in the form (\footnote{See H. Cramér, \textit{l.\,c.}, p.\ 17.})
\be
V(x)=\beta U_1(x)+(1-\beta)U_2(x)
\label{eqen26}
\ee
with $0\leq \beta \leq 1$ where $U_1(x)$ and $U_2(x)$ are two cumulative distribution functions such that \textit{almost everywhere}
\[
\begin{split}
U_1(x) &=\int_{-\infty}^x U'_1(y) dy,\\
U'_2(x) &=0.
\end{split}
\]

\switchcolumn*

\begin{conditionth}
\label{condfrB}
Dans la décomposition \eqref{eqfr26} de $V(x)$, on a $\beta>0$.
\end{conditionth}

\switchcolumn

\begin{conditionth}
\label{condenB}
In the decomposition \eqref{eqen26} of $V(x)$, we have $\beta>0$. \emph{\endnote{This condition requires the CDF $V(x)$ to have a smooth part with density $U'_1(x)$ and discrete part represented by $U_2(x)$.}}
\end{conditionth}

\switchcolumn*

Si $V(x)$ satisfait à la condition B, on voit immédiatement qu'il en est de même pour la fonction $\bV(x)$ définie par \eqref{eqfr6}. Pour la fonction $Q_n(y)$ définie par \eqref{eqfr22} on a alors (\footnote{Voir H. Cramér, \textit{l.\,c.}, p.\ 81. Il ne résulte pas immédiatement du théorème cité que la constante $K$ peut être prise indépendante de $h$. En parcourant la démonstration du théorème on s'assure cependant sans difficulté qu'il en est bien ainsi.})
\[
|Q_n(y)|< \frac{K}{\sqrt{n}}
\]
pour tout $n>1$, pour tout $y$ réel et pour tout $h$ de module suffisamment petit, la constante $K$ étant indépendante de $n$, $y$ et $h$.

\switchcolumn

If $V(x)$ satisfies Condition~\ref{condenB}, we see immediately that the same holds for the function $\bV(x)$ defined in \eqref{eqen6}. For the function $Q_n(y)$ defined by \eqref{eqen22}, we then have (\footnote{See H.\ Cramér, \textit{l.\,c.}, p.\ 81. It does not immediately follow from this theorem that the constant $K$ can be taken independent of $h$. Following the proof of this theorem, we nevertheless see without difficulty that this indeed the case.})
\[
|Q_n(y)|< \frac{K}{\sqrt{n}}
\]
for all $n>1$, for all real $y$, and for all $h$ of sufficiently small modulus, the constant $K$ being independent of $n$, $y$ and $h$.

\switchcolumn*

En introduisant ce résultat dans les calculs du chapitre III, on voit tout de suite que le facteur $\log n$, qui intervient dans les évaluations, peut être partout omis. De même la condition $x=o(\frac{\sqrt{n}}{\log n})$, dont le seul but est d'assurer la relation $\frac{x\log n}{\sqrt{n}}\ra 0$ peut être remplacée par $x=o(\sqrt{n})$. On a donc le théorème suivant.

\switchcolumn

By introducing this result into the calculations of Chapter III, we see directly that the factor $\log n$, which intervenes in the estimates, can be omitted everywhere. Likewise, the condition $x=o(\frac{\sqrt{n}}{\log n})$, whose only goal is to ensure the relation $\frac{x\log n}{\sqrt{n}}\ra 0$, can be replaced by $x=o(\sqrt{n})$. We thus have the following theorem.

\switchcolumn*

\begin{theoremfr}
\label{thmfr5}
Si les deux conditions~\ref{condfrA} et \ref{condfrB} sont satisfaites, on peut remplacer dans les théorèmes~\ref{thmfr1} et~\ref{thmfr4} la condition $x=o(\frac{\sqrt{n}}{\log n})$ par $x=o(\sqrt{n})$. On peut aussi omettre le facteur $\log n$ qui apparaît dans l'évaluation du reste dans les théorèmes~\ref{thmfr1}, \ref{thmfr2} et \ref{thmfr3}.
\end{theoremfr}

\switchcolumn

\begin{theoremen}
\label{thmen5}
If Conditions~\ref{condenA} and \ref{condenB} are satisfied, then we can replace in Theorems~\ref{thmen1} and \ref{thmen4} the condition $x=o(\frac{\sqrt{n}}{\log n})$ by $x=o(\sqrt{n})$. We can also omit the factor $\log n$ which appears in the evaluation of the rest in Theorems~\ref{thmen1}, \ref{thmen2} and \ref{thmen3}.
\end{theoremen}

\switchcolumn*

On peut cependant aller plus loin et considérer aussi les valeurs de $x$ qui sont du même ordre de grandeur que $\sqrt{n}$. Considérons en effet la condition~\ref{condfrA}, et désignons par $A_1$ et $-A_2$ les bornes supérieures et inférieures des valeurs de $h$ telles que l'intégrale \eqref{eqfr5} converge. $A_1$ et $A_2$ sont certainement des quantités positives, qui peuvent être finies ou non. En tenant compte des relations \eqref{eqfr14} et \eqref{eqfr15} on voit que, pour $-A_2<h<A_1$, la quantité $\bm$ définie par \eqref{eqfr7} est une fonction continue et toujours croissant de $h$, qui s'annule pour $h=0$. Les deux limites
\[
\lim_{h\ra A_1-0} \bm=\sigma C_1,\qquad \lim_{h\ra-A_2+0}\bm=-\sigma C_2,
\]
existent donc, $C_1$ et $C_2$ ayant des valeurs positives, finies ou non. Pour tout $c$ donné dans l'intervalle $-C_2< c< C_1$, l'équation 
\be
\bm=\sigma c
\label{eqfr27}
\ee
a une seule racine $h$ dans l'intervalle $-A_2<h<A_1$, dont le signe est le même que celui de $c$.

\switchcolumn

We can however go further and also consider values of $x$ that are of the same order as $\sqrt{n}$. Consider indeed Condition~\ref{condenA} and let us denote by $A_1$ and $-A_2$ the upper and lower bounds of the values of $h$ such that the integral \eqref{eqen5} converges. $A_1$ and $A_2$ are certainly positive quantities, which can be finite or not. Considering the relations \eqref{eqen14} and \eqref{eqen15}, we see that, for $-A_2<h<A_1$, the quantity  $\bm$ defined by \eqref{eqen7} is a continuous function, always growing with $h$, which vanishes for $h=0$. The two limits 
\[
\lim_{h\ra A_1-0} \bm=\sigma C_1,\qquad \lim_{h\ra-A_2+0}\bm=-\sigma C_2,
\]
thus exist, $C_1$ and $C_2$ having positive values, finite or not \mbox{\endnote{The two limits define the region of convergence of the generating function of $V(x)$.}}. For all given $c$ in the interval $-C_2< c< C_1$, the equation
\be
\bm=\sigma c
\label{eqen27}
\ee
has a unique root $h$ in the interval $-A_2<h<A_1$, whose sign is the same as that of $c$.

\switchcolumn*

Soit maintenant $h$ un nombre quelconque donné dans l'intervalle $-A_2  < h< A_1$, et considérons l'identité \eqref{eqfr21}, où $h$ entre comme paramètre. Les conditions~\ref{condfrA} et \ref{condfrB} étant satisfaites, nous avons pour la fonction $Q_n(y)$ définie par \eqref{eqfr22} le développement suivant (\footnote{Voir H. Cramér, \textit{l.\,c.}, p. 81.})
\[
Q_n(y)=\left(\frac{p_2(y)}{n^{\frac{1}{2}}}+\frac{p_5(y)}{n}+\cdots + \frac{p_{3k-1}(y)}{n^{\frac{k}{2}}}\right)e^{-\frac{y^2}{2}}+O\left(\frac{1}{n^{\frac{k+1}{2}}}\right),
\]
où $k$ est un entier arbitraire, tandis que les $p_\nu$ sont des polynômes dont le degré coïncide avec l'indice, et tels que les $p_{2\nu}$ sont des polynômes pairs, les $p_{2\nu-1}$ des polynômes impairs.

\switchcolumn

Now let $h$ be an arbitrary number given in the interval $-A_2< h< A_1$ and consider the identity \eqref{eqen21}, where $h$ enters as a parameter. Conditions~\ref{condenA} and \ref{condenB} being satisfied, we have for the function $Q_n(y)$ defined by \eqref{eqen22} the following expansion\\ (\footnote{See H. Cramér, \textit{l.\,c.}, p.\ 81.})
\[
Q_n(y)=\left(\frac{p_2(y)}{n^{\frac{1}{2}}}+\frac{p_5(y)}{n}+\cdots + \frac{p_{3k-1}(y)}{n^{\frac{k}{2}}}\right)e^{-\frac{y^2}{2}}+O\left(\frac{1}{n^{\frac{k+1}{2}}}\right),
\]
where $k$ is an arbitrary integer, whereas $p_\nu$ are polynomials whose degree coincides with the index, and such that $p_{2\nu}$ are even polynomials and $p_{2\nu-1}$ odd polynomials \endnote{The original article contains $p_{3k-1}$ instead of $p_{3k-1}(y)$. This is corrected here in both French and English versions.}.

\switchcolumn*

On en déduit, en refaisant les calculs du chapitre III,
\[
\int_0^\infty e^{-h\bsigma\sqrt{n} y}d\bF_n(y)=\frac{1}{\sqrt{n}}\left[b_0+\frac{b_1}{n}+\cdots+\frac{b_{k-1}}{n^{k-1}}+O\left(\frac{1}{n^k}\right)\right]
\]
pour tout entier positif $k$, les coefficients $b_\nu$ dépendant de $h$. On a d'ailleurs $b_0=\frac{1}{h\bsigma\sqrt{2\pi}}$. En introduisant dans \eqref{eqfr21}, on obtient donc
\be
\begin{split}
1- & F_n\left(\frac{\bm\sqrt{n}}{\sigma}\right)\\
&=\frac{1}{\sqrt{n}}e^{-(h\bm-\log R)n}\left[b_0+\frac{b_1}{n}+\cdots+\frac{b_{k-1}}{n^{k-1}}+O\left(\frac{1}{n^k}\right)\right].
\label{eqfr28}
\end{split}
\ee

\switchcolumn

We deduce, by re-doing the calculations of Chapter III,
\[
\int_0^\infty e^{-h\bsigma\sqrt{n} y}d\bF_n(y)=\frac{1}{\sqrt{n}}\left[b_0+\frac{b_1}{n}+\cdots+\frac{b_{k-1}}{n^{k-1}}+O\left(\frac{1}{n^k}\right)\right]
\]
for all positive integer $k$, the coefficients $b_\nu$ depending on $h$. We have in particular $b_0=\frac{1}{h\bsigma\sqrt{2\pi}}$. By substituting in \eqref{eqen21}, we then obtain
\be
\begin{split}
1  - & F_n\left(\frac{\bm\sqrt{n}}{\sigma}\right)\\
&=\frac{1}{\sqrt{n}}e^{-(h\bm-\log R)n}\left[b_0+\frac{b_1}{n}+\cdots+\frac{b_{k-1}}{n^{k-1}}+O\left(\frac{1}{n^k}\right)\right].
\label{eqen28}
\end{split}
\ee

\switchcolumn*

Soit maintenant $c$ un nombre donné tel que $0<c<C_1$, et prenons $h$ égal à la racine (unique) positive de l'équation \eqref{eqfr27}. En introduisant cette valeur dans \eqref{eqfr28} et en posant 
\be
\alpha=h\bm-\log R
\label{eqfr29}
\ee
(où l'on voit facilement que $\alpha$ est toujours positif), on a le théorème suivant.

\switchcolumn

Now let $c$ be a number such that $0<c<C_1$ and let us take $h$ equal to the (unique) positive root of Equation \eqref{eqen27}. Introducing this value in \eqref{eqen28} and taking
\be
\alpha=h\bm-\log R
\label{eqen29}
\ee
(where we easily see that $\alpha$ is always positive), we have the following theorem.

\switchcolumn*

\begin{theoremfr}
\label{thmfr6}
Si les deux conditions~\ref{condfrA} et \ref{condfrB} sont satisfaites, on peut trouver un nombre positif $C_1$ (fini ou non) tel que, pour tout $c$ dans l'intervalle $0<c<C_1$, on ait
\[
1-F_n(c\sqrt{n})=\frac{1}{\sqrt{n}}e^{-\alpha n}\left[b_0+\frac{b_1}{n}+\cdots+\frac{b_{k-1}}{n^{k-1}}+O\left(\frac{1}{n^k}\right)\right],
\]
où $\alpha$ est donné par \eqref{eqfr27} et \eqref{eqfr29}. Ici $k$ est un entier positif arbitraire, et les $b_\nu$ sont indépendants de $n$, mais dépendent de $c$. En particulier, on a toujours $b_0>0$.
\end{theoremfr}

\switchcolumn

\begin{theoremen}
\label{thmen6}
If the Conditions~\ref{condenA} and \ref{condenB} are satisfied, we can find a positive number $C_1$ (finite or not) such that, for all $c$ in the interval $0<c<C_1$, we have
\[
1-F_n(c\sqrt{n})=\frac{1}{\sqrt{n}}e^{-\alpha n}\left[b_0+\frac{b_1}{n}+\cdots+\frac{b_{k-1}}{n^{k-1}}+O\left(\frac{1}{n^k}\right)\right],
\]
where $\alpha$ is given by \eqref{eqen27} and \eqref{eqen29}. Here $k$ is an arbitrary positive integer and $b_\nu$ are independent of $n$, but depend on $c$. In particular, we always have $b_0>0$.
\end{theoremen}

\switchcolumn*

Il y a évidemment une expression analogue pour $F_n(-c\sqrt{n})$ où $-C_2<c<0$.

\switchcolumn

There is obviously an analoguous expression for $F_n(-c\sqrt{n})$ where $-C_2<c<0$.

\switchcolumn*

Ce théorème donne lieu à une remarque intéressante. Si les conditions~\ref{condfrA} et~\ref{condfrB} sont satisfaites, il suit du théorème~\ref{thmfr1} (avec les compléments apportés par le théorème~\ref{thmfr5}) que l'on a, pour $x\ra\infty$, $x=o(\sqrt{n})$,
\[
\log[1-F_n(x)]-\log[1-\Phi(x)]=x^2 o(1)=o(\log[1-\Phi(x)]),
\]
d'où
\[
\log[1-F_n(x)]\sim\log[1-\Phi(x)].
\]

\switchcolumn

This theorem gives rise to an interesting remark. If Conditions~\ref{condenA} and \ref{condenB} are satisfied, it follows from Theorem~\ref{thmen1} (with the complements given by Theorem~\ref{thmen5}) that we have, for $x\ra\infty$ and  $x=o(\sqrt{n})$,
\[
\log[1-F_n(x)]-\log[1-\Phi(x)]=x^2 o(1)=o(\log[1-\Phi(x)]),
\]
whence
\[
\log[1-F_n(x)]\sim\log[1-\Phi(x)].
\]

\switchcolumn*

D'autre part, pour $x=c\sqrt{n}$, on déduit du théorème~\ref{thmfr6}
\[
\log[1-F_n(x)]\sim-\frac{\alpha}{c^2}x^2\sim\frac{2\alpha}{c^2}\log[1-\Phi(x)],
\]
où, en général, la constante $\frac{2\alpha}{c^2}$ diffère de l'unité.

\switchcolumn

Moreover, for $x=c\sqrt{n}$, we deduce from Theorem~\ref{thmen6}
\[
\log[1-F_n(x)]\sim-\frac{\alpha}{c^2}x^2\sim\frac{2\alpha}{c^2}\log[1-\Phi(x)],
\]
where, in general, the constant $\frac{2\alpha}{c^2}$ is different from 1.

\switchcolumn*

\section*{Chapitre VI}

\switchcolumn

\section*{Chapter VI}

\switchcolumn*

Considérons maintenant une variable aléatoire $Z_t$, fonction d'un paramètre continu $t$, qu'on peut interpréter comme signifiant par exemple le temps. Supposons que l'accroissement $\Delta Z_t=Z_{t+\Delta t}-Z_t$ soit toujours une variable aléatoire indépendante de $Z_t$, et que la loi de répartition de $\Delta Z_t$ ne dépende ni de $t$ ni de $Z_t$, mais seulement de $\Delta t$. On dit alors que la variable $Z_t$ définit un \emph{processus stocastique homogène}. Supposons encore que la valeur moyenne $E(Z_t)$ s'annule pour tout $t$, et que $E(Z^2_t)$ soit toujours fini.

\switchcolumn

Consider now a random variable $Z_t$, a function of the continuous parameter $t$, which we can interpret, for example, as the time. Suppose that the increment \endnote{``Accroissement'' could be translated as ``difference'' or ``variation''. We use here the technical term ``increment''.} $\Delta Z_t=Z_{t+\Delta t}-Z_t$ is always a random variable independent of $Z_t$ and that the law  of $\Delta Z_t$ does not depend on $t$ nor on $Z_t$, but only on $\Delta t$. We then say that the variable $Z_t$ defines a \emph{homogeneous stochastic process}. Suppose furthermore that the expectation $E(Z_t)$ vanishes for all $t$ and that $E(Z^2_t)$ is always finite.

\switchcolumn*

Il résulte alors d'un théorème de M.\ Kolmogoroff (\footnote{Sulla forma generale di un processo stocastico omogeneo. \textit{Rend.~R.~Accad.~Lincei}, (6), 15 (1932), p.\ 805 et p.\ 866.}) qu'on peut assigner une constante $\sigma_0^2\geq 0$ et une fonction $S(x)$ bornée, jamais décroissante et continue au point $x=0$, avec les propriétés suivantes. Posons
\[
\begin{split}
\sigma_1^2 &= S(+\infty)-S(-\infty),\\
\sigma^2 &=\sigma_0^2+\sigma_1^2,\\
F(x,t)	 &= \textrm{Prob}(Z_t\leq \sigma x\sqrt{t}),\\
f(z,t)	&=\int_{-\infty}^\infty e^{izy}dF(y,t).
\end{split}
\]

\switchcolumn

It then follows from a theorem of Kolmogorov~\endnote{``Kolmogorov'' is used instead of ``Kolmogoroff''.} (\footnote{Sulla forma generale di un processo stocastico omogeneo.\ \textit{Rend.\ R.\ Accad.\ Lincei}, (6), 15 (1932), p.\ 805 and p.\ 866.}) that we can assign a constant $\sigma_0^2\geq 0$ and a bounded function $S(x)$, never decreasing and continuous at $x=0$, with the following properties. Define
\[
\begin{split}
\sigma_1^2 &= S(+\infty)-S(-\infty),\\
\sigma^2&=\sigma_0^2+\sigma_1^2,\\
F(x,t)&= \textrm{Prob}(Z_t\leq \sigma x\sqrt{t}),\\
f(z,t)&=\int_{-\infty}^\infty e^{izy}dF(y,t).
\end{split}
\]

\switchcolumn*

(Les différentielles devront toujours être prises par rapport à la variable $y$). Alors on a
\[
E(Z_t^2)=\sigma^2t
\]
et
\[
\log f(z,t)=-\frac{\sigma_0^2}{2\sigma^2}z^2+\frac{1}{\sigma^2}\int_{-\infty}^\infty\frac{e^{izy}-1-izy}{y^2}dS(\sigma y \sqrt{t}).
\]

\switchcolumn

(Derivatives must always be taken with respect to the variable $y$). We then have
\[
E(Z_t^2)=\sigma^2t
\]
and
\[
\log f(z,t)=-\frac{\sigma_0^2}{2\sigma^2}z^2+\frac{1}{\sigma^2}\int_{-\infty}^\infty\frac{e^{izy}-1-izy}{y^2}dS(\sigma y \sqrt{t}).
\]

\switchcolumn*

Il s'ensuit sans difficulté
\[
\lim_{t\ra\infty} f(z,t)=e^{-\frac{z^2}{2}},
\]
ce qui implique
\[
\lim_{t\ra\infty} F(x,t)=\Phi(x)
\]
pour tout $x$ réel indépendant de $t$.

\switchcolumn

There follows without difficulty
\[
\lim_{t\ra\infty} f(z,t)=e^{-\frac{z^2}{2}},
\]
which implies
\[
\lim_{t\ra\infty} F(x,t)=\Phi(x)
\]
for all real $x$ independent of $t$.

\switchcolumn*

Ici encore, on peut donc poser le problème d'étudier le comportement de $F(x,t)$ lorsque $x$ peut varier avec $t$, en tendant vers $+\infty$ ou vers $-\infty$ quand $t$ tend vers l'infini. Ce problème n'est en réalité qu'un cas particulier du problème dont nous nous sommes occupés dans les chapitres précédents.

\switchcolumn

Here again we can thus consider the problem of studying the behavior of $F(x,t)$ when $x$ can vary with $t$, going to $+\infty$ or to $-\infty$ when $t$ goes to infinity. This problem is nothing but a particular case of the problem that we have considered in the previous chapters.

\switchcolumn*

Supposons que l'intégrale
\be
\int_{-\infty}^\infty e^{hy}dS(y)
\label{eqfr30}
\ee
converge pour tout $h$ de module suffisamment petit, et posons
\[
\bS(x)=\int_{-\infty}^x e^{hy} dS(y)
\]
Alors il existe une variable $\bZ_t$ liée à un processus stocastique homogène, dont la répartition est définie au moyen de $\sigma_0^2$ et $\bS(x)$ de la même manière que la répartition de $Z_t$ a été définie par $\sigma_0^2$ et $S(x)$. Désignons par $\bsigma_1^2$, $\bsigma^2$, $\bF(x,t)$ et $\bar{f}(z,t)$ les quantités analogues aux précédentes formées en partant de $\sigma_0^2$ et $\bS(x)$. Définissons ici les quantités $\bm$ et $R$ en posant
\[
\bm=\sigma_0^2 h+\int_{-\infty}^\infty \frac{e^{hy}-1}{y}dS(y),
\]
\[
\log R=\frac{1}{2}\sigma_0^2 h^2+\int_{-\infty}^\infty \frac{e^{hy}-1-hy}{y^2}dS(y).
\]

\switchcolumn

Suppose that the integral
\be
\int_{-\infty}^\infty e^{hy}dS(y)
\label{eqen30}
\ee
converges for all $h$ of sufficiently small modulus, and let
\[
\bS(x)=\int_{-\infty}^x e^{hy} dS(y).
\]
Then there exists a variable $\bZ_t$ linked to a homogeneous stochastic process whose distribution function is defined by means of  $\sigma_0^2$ and $\bS(x)$, in the same way that the distribution of $Z_t$ was defined by $\sigma_0^2$ and $S(x)$. Denote by $\bsigma_1^2$, $\bsigma^2$, $\bF(x,t)$ and $\bar{f}(z,t)$ the quantities similar to those formed before by starting from $\sigma_0^2$ and $\bS(x)$. Define here the quantities $\bm$ and $R$ as 
\[
\bm=\sigma_0^2 h+\int_{-\infty}^\infty \frac{e^{hy}-1}{y}dS(y),
\]
\[
\log R=\frac{1}{2}\sigma_0^2 h^2+\int_{-\infty}^\infty \frac{e^{hy}-1-hy}{y^2}dS(y).
\]

\switchcolumn*

On démontre alors par un calcul analogue à celui du chapitre II l'identité suivante qui a lieu pour tout $h$ réel appartenant au domaine de convergence de l'intégrale \eqref{eqfr30}\,:
\[
1-F\left(\frac{\bm\sqrt{t}}{\sigma},t\right)=R^te^{-h\bm t}\int_0^\infty e^{-h\bsigma\sqrt{t}y}d\bF(y,t).
\]
Cette identité est, comme on le voit, tout à fait analogue à l'identité \eqref{eqfr21}. On peut aussi s'en servir d'une manière absolument analogue, en démontrant des théorèmes sur la fonction $F(x,t)$ qui sont parfaitement analogues aux théorèmes~\ref{thmfr1}-\ref{thmfr6} sur la fonction $F_n(x)$. La seule différence est que le paramètre discontinu $n$ a été remplacé par le paramètre continu $t$. Dans les conditions~\ref{condfrA} et \ref{condfrB}, on doit remplacer la fonction $V(x)$ par la fonction $S(x)$ considérée dans ce chapitre (\footnote{La condition \ref{condfrB} peut être remplacée par une autre condition plus générale que celle obtenue de la manière indiquée. Voir H. Cramér, \textit{l.\,c.}, p. 99. -- Un théorème contenu dans notre théorème~\ref{thmfr6} a été énoncé, pour un cas particulier important du processus homogène, par F. Lundberg, Försäkringsteknisk riskutjämning, Stockholm 1926-1928, et démontré par F. Esscher, \textit{l.\,c.}}).

\switchcolumn

We then demonstrate through a calculation similar to that of Chapter II the following identity which holds for all real $h$ in the domain of convergence of the integral \eqref{eqen30}:
\[
1-F\left(\frac{\bm\sqrt{t}}{\sigma},t\right)=R^t e^{-h\bm t}\int_0^\infty e^{-h\bsigma\sqrt{t}y}d\bF(y,t).
\]
This identity is, as we see, totally analoguous to the identity \eqref{eqen21}. We can also use it in a manner absolutely analoguous, by proving the theorems on the function $F(x,t)$ which are perfectly analoguous to Theorems~\ref{thmen1}-\ref{thmen6} about the function $F_n(x)$. The only difference is that the discontinuous parameter $n$ has been replaced by the continuous parameter $t$. In the Conditions~\ref{condenA} and \ref{condenB}, we must replace the function $V(x)$ by the function $S(x)$ considered in this chapter (\footnote{Condition~\ref{condenB} can be replaced by another condition that is more general than the one indicated. See H.\ Cramér, \textit{l.\,c.}, p.\ 99. -- A theorem contained in our Theorem~\ref{thmen6} has been stated, for a particular important homogeneous case, by F.\ Lundberg, Försäkringsteknisk riskutjämning, Stockholm 1926-1928, and demonstrated by F.\ Esscher, \textit{l.\,c.}}).

\switchcolumn*

\end{paracol}


\newgeometry{margin=1.25in,top=1.5in,bottom=1in}
\fontsize{11pt}{12.5pt}\selectfont


\theendnotes

\vfill

{\footnotesize\noindent This document was compiled with \LaTeX\ using the \texttt{paracol} package for aligning the French and English texts, after a custom Perl script meshed them from separate files.}

\end{document}